\documentclass[11pt]{article}
\usepackage{amsmath, amsthm, amsfonts, amssymb}
\usepackage{color}

\hoffset= - 0.5 cm \voffset=  0.1 cm \numberwithin{equation}{section}

\newtheorem{theorem}{Theorem}[section]
\newtheorem{lemma}[theorem]{Lemma}
\newtheorem{proposition}[theorem]{Proposition}
\newtheorem{corollary}[theorem]{Corollary}
\newtheorem{remark}[theorem]{Remark}
\newtheorem{definition}[theorem]{Definition}
\newtheorem{example}[theorem]{Example}
\newtheorem{hypothesis}[theorem]{Hypothesis}

\newcommand{\bth}{\begin{theorem}}
\renewcommand{\eth}{\end{theorem}}
\newcommand{\bpr}{\begin{proposition}}
\newcommand{\epr}{\end{proposition}}
\newcommand{\bco}{\begin{corollary}}
\newcommand{\eco}{\end{corollary}}
\newcommand{\ble}{\begin{lemma}}
\newcommand{\ele}{\end{lemma}}

\newcommand{\bpf}{\begin{proof}}
\newcommand{\epf}{\end{proof}}
\newcommand{\bex}{\begin{example}}
\newcommand{\eex}{\end{example}}
\newcommand{\bdf}{\begin{definition}}
\newcommand{\edf}{\end{definition}}
\newcommand{\bre}{\begin{remark}}
\newcommand{\ere}{\end{remark}}

\newcommand{\beq}{\begin{equation}}
\newcommand{\eeq}{\end{equation}}
\newcommand{\bal}{\begin{aligned}}
\newcommand{\eal}{\end{aligned}}
\newcommand{\beqr}{\begin{eqnarray*}}
\newcommand{\eeqr}{\end{eqnarray*}}

\def\({\Bigl(}\def\){\Bigr)}
\def\P{{\mathbb P}}
\def\R{{\mathbb R}}

\def\E{{\mathbb E }}
\def\N{\mathbb N}

\def\lan{{\langle}}
\def\ran{{\rangle}}
\def\hh{{\vskip 2mm \noindent }}
 \def\vv{\vskip 1mm \noindent}

\def\ti{\tilde}

\begin{document}

\title {\Huge Structural properties of   semilinear  SPDEs  driven by
 cylindrical  stable processes
\footnote{\sl
 The paper is almost identical
with the paper published under the same title in Probab. Theory
Related Fields  (see \cite{PZ7}) with the exception of some
constants in particular in  Theorem 4.16 and Hypothesis 5.6. We also
thank Lihu Xu for indicating an error in our previous calculations.}
 }

\date{}

\maketitle


\begin{center}
 Enrico Priola \footnote
 { \noindent \! \! \!
  Supported by the M.I.U.R. research projects Prin 2004 and
 2006 ``Kolmogorov equations'' and
 by the Polish Ministry of Science and Education project 1PO 3A 034
29 ``Stochastic evolution equations with L\'evy noise''.}

\vspace{ 2 mm } {\small  \it Dipartimento di Matematica,
Universit\`a di Torino, \par
  via Carlo Alberto 10,  \ 10123, \ Torino, Italy. \par
 e-mail \ \ enrico.priola@unito.it }
\\ \ \\
Jerzy Zabczyk   \footnote { \noindent \! \! Supported by the
Polish Ministry of Science and Education project 1PO 3A 034 29
``Stochastic evolution equations with \! \! L\'evy noise''.}

\vspace{ 2 mm} {\small  \it Instytut Matematyczny,   Polskiej
Akademii Nauk,
\par ul. Sniadeckich 8, 00-950, \ Warszawa, Poland.\par
 e-mail \ \ zabczyk@impan.gov.pl }
\end{center}


\vskip 7mm
 \noindent {\bf Mathematics  Subject Classification (2000):} \
   60H15,
  60J75, 47D07, 35R60.
 \ \par \ \par

 \noindent {\bf Key words:} Stochastic PDEs with jumps,
 strong Feller property, regularity of trajectories.

\vspace{1.5 mm}

\noindent {\bf Abstract:} We consider a class of semilinear
 stochastic evolution equations   driven by an additive  cylindrical stable
 noise.
    We investigate  structural properties of the solutions
 like  Markov, irreducibility, stochastic continuity,  Feller and
 strong Feller properties,  and
 study   integrability  of  trajectories. The obtained
 results  can be applied to
 semilinear stochastic heat equations with Dirichlet boundary
 conditions and
  bounded and Lipschitz
  nonlinearities.

\section{Introduction }

 The paper is concerned with  structural properties   of
 solutions to nonlinear stochastic equations
 \beq
\label{ab}
 dX_t  = A X_t dt + F(X_t)dt +
    dZ_t, \;\;\;\;  t \ge 0, \;\; \; X_0  = x \in  H,
 \eeq
   in a real separable Hilbert space  $H$ driven by
    an infinite dimensional stable process $Z= (Z_t)$. In
    particular, we study
     Markov, irreducibility, stochastic continuity, Feller and
   strong Feller properties for the solutions, and
 investigate  integrability  of  trajectories.
  The main results are  gradient estimates for the associated
   transition semigroup (see Theorem \ref{grad} when $F=0$
   and Theorem \ref{main} in the general case),
 from which we deduce the strong Feller property,
  and a theorem
  on  time regularity of trajectories (see Theorem \ref{kw}).

  To cover interesting cases, we consider processes  $Z$ which
    take values in a Hilbert space $U$ usually
greater than $H$.
    Moreover  $A : $ dom$(A) \subset
   H \to H $ is a linear possibly unbounded operator
  which generates a $C_0$-semigroup $(e^{tA})$ on
   $H$
     and $F : H \to H$
    denotes a Lipschitz continuous and bounded function.

  In the  case  when $Z$ is a Wiener process the theory of
 equations \eqref{ab} is well understood.
  The situation changes completely
  in  the stable noise case and new phenomena appear.
  For instance, even in the linear case
 $F=0$, it is not known when solutions of \eqref{ab}
  have c\`adl\`ag trajectories.
 That  lack of c\`adl\`ag regularity   is possible was noted in
(\cite[Proposition 9.4.4]{PeZ}) in a similar situation.
  Another
difficulty is related to the fact that general necessary and
sufficient conditions for absolute continuity of stable measures
on Hilbert spaces  do not exist.

  We restrict our considerations to
 SPDEs with additive noise
as even in this case some new phenomena, related to the
cylindrical L\'evy noise, appear.{ We hope that the results
presented here will
 form a proper starting
point to treat general equations with  multiplicative L\'evy
perturbations}.

   In this paper
   we consider
  a {\it cylindrical $\alpha$-stable process}  $Z= (Z_t)$,  $\alpha \in
 (0,2)$,
   defined  by
the orthogonal expansion
\begin{equation}\label{ee}
 Z_t = \sum_{n \ge 1} \beta_n Z_t^n e_n,\;\;\;\; t \ge 0,
\end{equation}
 where  $(e_n)$ is an orthonormal basis of $H$ and
    $Z_t^n$ are independent,  real
 valued, normalized, symmetric $\alpha$-stable
  processes defined on a fixed stochastic basis. Moreover,
    $(\beta_n)$ is a given, possibly unbounded,
   sequence of {\it positive} numbers.

 The results of the paper
 apply  to stochastic heat
 equations with Dirichlet boundary conditions
 \beq \label{heat}
  \left\{
\bal dX(t, \xi)  & = \big( \triangle X(t, \xi) +
 f(X(t, \xi)) \big ) \,dt  + \,dZ(t, \xi), \;\;\; t>0,\\
      X(0, \xi) & = x(\xi), \;\;\; \xi \in D,
\\ X(t, \xi) &=0, \;\; t>0, \;\;\; \xi \in \partial D,
\eal \right. \eeq
 in a given bounded domain $D \subset \R^d$
   having Lipschitz-continuous
   boundary $\partial D$. Here $x(\xi) \in H= L^2(D)$, $f:
   \R \to \R$ is bounded and Lispchitz continuous and
 the noise $Z $ is a cylindrical
    $\alpha$-stable process  of the form
 \eqref{ee}, where
 $(e_n)$ is a basis of eigenfunctions for the Laplace
  operator  $\triangle$
  (with  Dirichlet boundary conditions).

\smallskip  Irreducibility and strong Feller property can
 be  used to establish uniqueness of an invariant measure for the
  solutions of \eqref{ab} through  a well known  approach
   based on the so called
   Doob theorem
   (see \cite[Theorem 4.2.1]{DZ1}). However our results indicate
   that solutions of \eqref{ab} with non-Gaussian
    noise are less regular than those with a  Wiener process.
  Thus to cover equations with L\'evy noise, having only
   a finite number of
 modes or with modes vanishing rapidly (compare with \cite{HM1})
  one needs an extension of the methods
  developed in \cite{HM1}, \cite{HM2} and \cite{HM3}.

\vskip 1mm  After short  Preliminaries, concerned with notations
and basic definitions, in Section 3, we deal   with real and
Hilbert space
 valued  $\alpha$-stable random variables.
 The
most important result here is a
 necessary and sufficient
 condition for the
  absolute continuity of shifts of infinite products
 of symmetric $\alpha$-stable, one dimensional distributions
  (see Theorem
\ref{gh}).
 It is an improvement of an old result by Zinn (see
\cite{Zi}) with a direct proof.  Section 4 is concerned with
linear equations \beq \label{li} dX_t= A X_t dt +  dZ_t, \;\; \; t
\ge 0, \;\;\;
 X_0=x \in H
 \eeq
  (see also  \cite{Ch},
   \cite{BRS},  \cite{FR}, \cite{Dawson}, \cite{RW} and
   \cite{PeZ}).  We  assume
    that vectors $(e_n)$
from the representation \eqref{ee} are eigenvectors of $A$.

In Proposition \ref{bo1}  we give if and only if
  conditions under
 which $X$, the solution of (\ref{li}),
  takes values in $H$,
  and establish its  measurability and markovianity.
Then we deal with the time regularity of trajectories. The main
result here is Theorem \ref{kw}, which establishes
 stochastic continuity of the solution and  integrability of
 its trajectories. Better regularity, like right or left continuity of
trajectories is established here in   very special cases  and {
 is an open question for general equations}. Note that
  in \cite{FR} it is proved that trajectories of $(X_t^x)$
    are
  c\`adl\`ag only in some enlarged Hilbert space $U$ containing
  $H$. This lack of time regularity
introduces additional difficulties into the theory
 (see also Proposition
 9.4.4 in \cite{PeZ}). We establish  also irreducibility of the
 solution. Theorem \ref{ass} gives conditions under which all
transition laws of $X$ are equivalent and establishes a formula
for the densities. Moreover
           (see Theorem \ref{grad})  under the
          assumptions of
Theorem \ref{ass}, the transition semigroup corresponding to $X$ is
not only strong Feller but transforms bounded measurable functions
onto Fr\'echet differentiable functions with continuous derivative.
Important gradient estimates are established as well. Theorems on
Ornstein-Uhlenbeck processes are based on results
 about
stable measures established
 in Section 3.

  Section 5  is devoted to nonlinear equations \eqref{ab}.
    Markov property and
 irreducibility require special
  attention due to the lack of c\`adl\`ag
 regularity of the trajectories.
 They are  established
  in Theorems \ref{exist} and \ref{irre}. Then
estimates of Section 4 are  used to establish the strong Feller
property of the solution to the nonlinear equation
 (see Theorem \ref{main}).   Here we assume  $\alpha \in (1,2)$. The main tool  is the
 so called mild
  version of the Kolmogorov equation and Galerkin's approximation.
  It is proper to add that the classical approach to get
  strong Feller using
 the  Bismut-Elworthy-Li formula
 is not available in the non-Gaussian
case.
 A related formula, but requiring a non trivial Gaussian
 component in the L\'evy noise,
 was established in finite dimensions in \cite{PZ6}.

\section {Preliminaries}

   $H$ will denote  a real separable
   Hilbert with inner product
    $\langle \cdot, \cdot \rangle$ and norm
$|\cdot|$. By ${\cal L} (H)$ we denote the Banach space
 of all bounded
 linear operators from $H$ into $H$ endowed with the operator norm
  $\|  \cdot \|_{{\cal L}(H)}$.
 We will fix an  orthonormal basis
   $(e_n)$ in $H$. Through
    the basis $(e_n)$ we will often identify $H$
 with  $l^2$.
  More generally, for a given sequence $\rho= (\rho_n)$
 of real numbers, we set
   \beq \label{rho}
 l^2_{\rho} = \{ (x_n) \in \R^{\infty} \,\, :\, \sum_{n \ge 1} x_n^2
 \rho_n^2 < \infty \},
\eeq
 where $\R^{\infty} = \R^{\N}$.
  $l^2_{\rho}$ becomes a separable Hilbert space with the inner
 product: $\langle x,y\rangle = \sum_{n \ge 1} x_n y_n
 \, \rho_n^2$, for $x = (x_n)$, $y = (y_n) \in l^2_{\rho}$.

\smallskip The space  ${   C}_b ^{}(H )$ (resp. ${ B}_b(H)$)
 stands for the
 Banach space of all  real, continuous (resp. Borel) and bounded  functions \ $f:
 H \to \R $, endowed with the supremum norm:
 \  $
 \| f \|_0 \, =  \sup _{x \in H } | f(x) |.$

{\vskip 1mm} The space ${  C}_b^{k }(H)$, $k \ge 1$,
  is the set of
  all $k$-times differentiable functions $f: H \to \R$,
   whose Fr\'echet  derivatives
 $D^i f$, $1 \le i \le k,$   are  continuous and bounded on $H$,
 up to the order $k $.
   Moreover we set ${ C}_b^{\infty }(H)$ $= \cap_{k \ge 1}$
 ${ C}_b^{k }(H).$

{\vskip 1mm} Let us recall that  a  L\'evy process $(Z_t)$ with
values in $H$ is an
 $H$-valued process defined  on
 some stochastic basis
  $(\Omega, {\cal F}, ({\cal F}_t)_{t \ge 0}, \P)$,
 having stationary independent increments,
  c\`adl\`ag  trajectories,
  and such that $Z_0 =0$, $\P$-a.s..

One has that
 \begin{equation}\label{zt}
\E [e^{i \langle Z_t , s \rangle }]  = \exp ( - t \psi (s) ), \; s
\in H,
 \end{equation}
where the exponent  $\psi$ can be expressed by the following
infinite dimensional  L\'evy-Khintchine formula,
 \begin{equation}\label{psi} \psi(s)= \frac{1}{2} \langle
Q s,s \rangle - i \langle a, s\rangle - \int_{H} \Big(  e^{i \langle
s,y \rangle }  - 1 - \frac { i \langle s,y \rangle} {1 + |y|^2 }
\Big ) \nu (dy), \;\; s \in H.
 \end{equation}
 Here $Q$ is a  symmetric non-negative trace class operator on
$H$, $a \in H$ and $\nu$ is the  L\'evy measure or the jump
intensity measure  associated to $(Z_t)$
 (see  \cite{sato} and \cite{PeZ}).

\vskip 1mm    According to Proposition \ref{ll}
 (see also Remark \ref{levy})
  our cylindrical $\alpha$-stable process
 $Z$ appearing in \eqref{ee}
  is a L\'evy process
  taking values in the Hilbert space
   $U =l^2_{\rho}$, see \eqref{rho}, with a properly
 chosen weight ${\rho}$.

\section{Stable measures on Hilbert spaces}

 Here we gather and strengthen results on stable distributions
needed in the sequel.

\subsection{Stable  densities}

 Let us consider a  one
 dimensional, normalized,  symmetric $\alpha$-stable
distribution
  $\mu_{\alpha}$, $\alpha \in
 ]0,2]$, having  characteristic function
  \beq \label{al}
 \hat \mu_{\alpha} (s) = e^{- |s|^{\alpha}},\;\;\; s \in \R.
 \eeq
The
 {\it density of $\mu_{\alpha} $}, with respect to
the Lebesgue
 measure, will be denoted by
  $
 p_{\alpha}$.

\noindent We need to know
 the {\it precise asymptotic behaviour }
of the density $p_{\alpha}$, $\alpha \in (0,2)$.

We have that for any $\alpha \in (0,2)$, there exists
$C_{\alpha}>0$ such that
 \beq \label{as}
 p_{\alpha}(x) \, \sim \, \frac{C_{\alpha}}{x^{\alpha +1}},\;\;\;\;
 \mbox{as} \;\; x
 \to \infty,
\eeq
 see \cite{Z},  \cite[page 88]{sato}and
  \cite[pages 582-583]{Fe}.
 We need two lemmas about $p_{\alpha}$.  The first one  is
  straightforward.

 \ble \label{pp} Let $p_{\alpha}$ be the
  density of the  one dimensional
   $\alpha$-stable measure $\mu_{\alpha}$ in \eqref{al}.
    Then, for any $\alpha \in (0,2)$,
      $p_{\alpha} \in C^{\infty}(\R) \cap C_0(\R)$ and moreover
       (with $p_{\alpha}'(x) = \frac{d p_{\alpha}}{dx}$)
 \beq \label{as1}
 x^2 p_{\alpha}'(x) \in L^{\infty}(\R).
 \eeq
\ele
 \bpf Let $p = p_{\alpha}$.
 It is well known that $p \in C^{\infty}(\R) \cap
C_0(\R)$
 (see, for instance, \cite[Chapter 1]{sato}).
 To get the second assertion, we use the inversion Fourier formula and
  integrate by parts,
 \beqr
 x^2 p'(x ) = -\frac{i}{2 \pi}\int_{\R} x^2 e^{-ixy} y e^{-|y|^{\alpha}} dy
\\
= \frac{x}{2 \pi}  \int_{\R} \frac{d}{dy} \big ( e^{-ixy}\big) y
e^{-|y|^{\alpha}} dy
 =  \frac{x}{2 \pi} \int_{\R} e^{-ixy}
 e^{-|y|^{\alpha}} \Big (  \alpha |y|^{\alpha} -1
\Big) dy
\\
= \frac{-i  \alpha^2}{2 \pi}  \int_{\R} e^{-ixy}
 e^{-|y|^{\alpha}}  \frac{y} {|y|^{2 - \alpha}}
 dy
 \, + \, \frac{i \alpha}{2 \pi} \int_{\R} e^{-ixy}
 e^{-|y|^{\alpha}} y \Big(
\frac{ \alpha } {|y|^{2 - 2 \alpha}}
   -\frac{1} {|y|^{2 - \alpha}}
   \Big)dy.
\eeqr
 From this formula the assertion is clear.
\epf


\ble \label{gj} Let us consider  the function
$$
g(x) =  1-\int_{\R} p^{1/2}_{\alpha}(z ) \, p^{1/2}_{\alpha}\big(z
- x \big)\,dz, \;\;
 \;\;\; x \in (-1,1).
$$
 We have
 \beq \label{g}
 g(x) \, \sim \, c_{\alpha} x^2\;\;\; as \;\; x \to 0, \;\;\;
  \mbox{where} \;\;\; c_{\alpha} =  \frac{1}{8}
 \int_{\R}  \frac{p_{\alpha}'(z)^2}{p_{\alpha}(z )}
 \,dz,\;\;\; \alpha \in (0,2).
  \eeq
\ele \bpf Let $p = p_{\alpha}$. Clearly $g(0)=0$.
  In order to prove \eqref{g}
  we will apply Hopital's rule. To this purpose we prove that
   $g$ is twice differentiable, with $g'(0)=0$ and $g''(0) \not
   =0$.
We have, for $|x|<1$,
$$
 g'(x) =  \frac{1}{2}
 \int_{\R}  p^{1/2}(z ) \, \frac{1}{
 p^{1/2}\big(z - x \big)} \, p'(z-x)\,dz.
$$
 The differentiation is justified by \eqref{as} and \eqref{as1},
  using also the fact that $p$ is a positive function on $\R$.
 We only point out  the following useful estimate:   for any $M>1$, there
    exists $c>0$, such that,
    for any $x \in (-1,1)$, $|z|>M$,
$$
p^{1/2}(z) \frac{1}{
 p^{1/2}\big(z - x \big)} \, |p'(z-x)| \le \frac{c}{|z|^{1/2 + \alpha /2}}
 \frac{(|z|+1)^{1/2 + \alpha/2}}{|z -1|^2}.
$$
 We also get  $g'(0)
  =  \frac{1}{2}
 \int_{\R}  p'(z)\,dz =0$. We show now that there exists the
 second derivative of $g$. To this purpose, we write
$$
g'(x)=  \frac{1}{2}
 \int_{\R}  p^{1/2}(z + x) \, \frac{1}{
 p^{1/2}\big( z \big)} \, p'(z)\,dz.
$$
 We have, for any $x \in (-1,1)$,
$$
g''(x)=  \frac{1}{4}
 \int_{\R}  \frac{p'(z+x)}{p^{1/2}(z + x)} \, \frac{1}{
 p^{1/2}\big( z \big)} \, p'(z)\,dz.
$$
 The differentiation can be done, since, for any $M>1$,
 there exists $c'>0$, such that,
    for any $x \in (-1,1)$, $|z|>M$,
$$
\frac{|p'(z+x)|}{p^{1/2}(z + x)} \, \frac{|p'(z)|}{
 p^{1/2}\big( z \big)} \,  \le \,
 \frac{c' \, (|z|+1)^{1 + \alpha}}{|z|^2 \, |z -1|^2}.
$$
We have also that
$$
g''(0) =  \frac{1}{4}
 \int_{\R}  \frac{p'(z)^2}{p^{}(z )}
 \,dz .
$$
 and so \eqref{g} is proved.
\epf

\subsection{Supports of stable measures }

 Let us consider a sequence $(\xi_n)$ of
  independent real random variables, having
 the same law $\mu_{\alpha}$ and defined on a probability space
$(\Omega, {\cal F}, \P)$.
 Take nonnegative numbers $q_n$ and consider
the random variable
 \beq \label{cio}
 \xi = (q_1 \xi_1, \ldots, q_n \xi_n, \ldots )
 \eeq
with values in $\R^{\infty}$. We start with a preliminary result,
which is a special case of \cite[Corollary 2.4.2]{Kw}.
 We provide a proof for the sake of completeness.

\bpr \label{ll} For any $\alpha \in ]0,2]$,
 the random variable
 $\xi $ in \eqref{cio} takes values in $l^2$, $\P$-a.s., if and
 only if
 \beq \label{if}
 \sum_{n \ge 1} q_n^{\alpha} < \infty.
 \eeq
  If, in addition to \eqref{if}, $q_n >0$, $n =1, 2, \ldots, $
 then the support of the law of $\xi$ is $l^2$.
\epr \bpf
 We will use the
   following  theorem (see, for instance \cite{K},
   page 70-71): { let $U_n$
   be a sequence of independent and symmetric real random
   variables; then the following statements are equivalent:
    $\sum_{n \ge 1} U_n$ converge in
   distribution; $\sum_{n \ge 1} U_n$
  converges $\P$-a.s.; $\sum_{n \ge 1} U_n^2 $ converges
  $\P$-a.s..}

\vv We have:
  $ \E  [e^{i \sum_{n=1}^N  q_n \xi_n h } ] = \prod_{n=1}^N
\E  [e^{i  q_n \xi_n h} ]$ $ = e^{-  \sum_{n=1}^N \, q_n^{\alpha}
\, |h|^{\alpha}},$  for any $N \in \N$, $h  \in \R$,

 Then it is clear that $ \sum_{k=1}^N  q_k \xi_k $ converges in
 distribution if and only if (\ref{if}) holds.
Moreover if (\ref{if})  holds,  then we have  convergence in
distribution to
  the random variable $\xi_1 \big (\sum_{k=1}^{\infty} \,
   q_k^{\alpha} \, )^{1/\alpha}  $.
    It follows that the series $\sum_{k \ge 1} q_k \xi_k$
   converges, $\P$-a.s., and also that
   \beq \label{vi}
 \sum_{k \ge 1}  q_k^2 \, \xi_k^2 < \infty,\;\;\; \; \;\;
 \P-a.s.,
\eeq
 and this proves the first part.

 To prove the second assertion, we fix an arbitrary ball $B
\subset l^2$, $B = B(y, r)$ with center
 in $y= (y_k) \in l^2$ and radius $r>0$. Using independence, we find
 \beqr
 \P \Big(\sum_{k \ge 1} (q_k\xi_k - y_k)^2 < r^2 \Big)
  \ge \P \Big(\sum_{k = 1}^N (q_k \xi_k - y_k)^2 < \epsilon \Big)
 \, \P \Big(\sum_{k > N} (q_k \xi_k - y_k)^2 < r^2 - \epsilon \Big).
   \eeqr
Now we use that the one dimensional measure $\mu_{\alpha}$ has a
positive density on $\R$. This implies that,  for any $N \in \N$,
$\epsilon
>0$,
$ \P \Big(\sum_{k = 1}^N (q_k \xi_k - y_k)^2 < \epsilon \Big) >0.
$

Since $\P (
 \, \sum_{k > N} (q_k\xi_k - y_k)^2 < r^2 - \epsilon) \to 1$,
 as $N \to \infty$, the assertion follows.
\epf

\subsection{Equivalence of shifts of stable measures}

 Here we   give necessary  and sufficient
  conditions in
  order that  shifts of infinite products of  one dimensional
  $\alpha$-stable distributions, are equivalent.
  Our theorem
  on equivalence
  strengthen an absolute continuity result of Zinn \cite{Zi}
 (see Remark \ref{zinn}) with a different proof
  which requires   Lemma \ref{gj}..

\bth \label{gh} Let us consider the $l^2$-random variable  $\xi $
in \eqref{cio} under the condition $q_k >0$, $k \ge 1$, and
$\sum_{k \ge 1} q_k^{\alpha} < \infty$.  Take arbitrary $u, v \in
l^2$ such that
 \beq \label{new}
  \sum_{k \ge 1} \frac{|u_k - v_k|^2}{q_k^2} <
\infty.
 \eeq  Then the law of the random variable
 $\xi + u$ and the one of $\xi + v$ are equivalent.

 In addition, if $\mu$ and $\nu$ denote the laws of
  $\xi + u $ and $\xi + v$ respectively,
 the density $\frac{d \mu }{d \nu }$ of
 $\mu $ with respect to $\nu$ is given by
 $$
\frac{d \mu }{d \nu} = \lim _{N \to \infty} \prod_{k=1}^N
  \frac{p_{\alpha} \Big(\frac{ z_k - u_k}{q_k}
 \Big)}{p_{\alpha} \Big(\frac{ z_k - v_k  }{q_k}
 \Big)} \;\;\;\;\;  in \;\; L^1 (\nu),\;\; \alpha \in ]0,2].
$$
 \eth
\bpf
   The result is well known  when $\alpha =2$ so
   let $p_{\alpha} = p$ with $\alpha \in (0,2)$.
   The measures  $\mu $ and $\nu$ can be seen as Borel
    product measures
in $\R^{\infty}$, i.e.,
$$
 \bal
& \mu = \prod_{k \ge 1} \mu^k,\;\;\;\;\;
   \nu = \prod_{k \ge 1} \nu^k, \;\;\; \mbox{where}
  \;\;  \mu^k, \; \nu^k \;\;  \mbox{have densities, respectively, }
\\
&  \frac{1}{q_k } p \Big(\frac{ z_k - u_k} {q_k}
 \Big)  \;\; \mbox{and} \;\;
  \frac{1}{q_k } p \Big(\frac{ z_k - v_k} {q_k}
 \Big).
 \eal
$$
  According to
\cite[Proposition 2.19]{DZ}, $\mu$ and $\nu$ are equivalent if and
only if
$$
 H(\mu, \nu) = \prod_{k \ge 1} \int_{\R} \Big (
  \frac{d \mu^k }{d \nu^k} \Big)^{1/2}
  \nu^k ( d  z_k) \,
 = \, \prod_{k \ge 1} \int_{\R} \Big (
  \frac{d \mu^k }{d z_k} \Big)^{1/2}
  \, \Big (
  \frac{d \nu^k }{d z_k} \Big)^{1/2}
   ( d  z_k)
\;  >0.
$$
Define, for any $k \ge 1,$
$$
\aligned a_k&=\int_{\R
}\biggl(\frac{d\mu^k}{dz_k}(z_k)\biggr)^{1/2}
\biggl(\frac{d\nu^k}{dz_k}(z_k)\biggr)^{1/2}dz_k
\\
&=\int_{\R}\biggl[p^{1/2} \Bigl(z_k-\frac{u_k}{q_k} \Bigr) p^{1/2}
\( z_k- \frac{v_k}{q_k}\)\biggr]\,dz_k \in (0, 1].
\endaligned
$$
 Note that
$$
\prod_{k \ge 1} a_k =\prod_{k \ge 1} (1-(1-a_k))=e^{\sum_{k \ge 1}
\, \ln(1-(1-a_k))}
$$
Note that, if $0 \le 1 - a \le  1/2$, then
$$
\ln(1-(1-a)) \ge (- 2 \log 2) \, (1-a).
$$
Consequently, if we prove that there exists $k_0$ such that, for
all $k \ge k_0$,  $0 \le 1-a_k \le 1/2$, then we get
$$
 \;\; \prod_{k \ge k_0}
a_k\ge e^{- 2 \log 2 \, \sum_{k \ge k_0} \, (1-a_k)}.
$$
Let us write
$$
 1-a_k=1-\int_{\R} p^{1/2}(z_k- \frac{u_k}{q_k} )p^{1/2}(z_k-
\frac{v_k}{q_k})\,dz_k = 1-\int_{\R} p^{1/2}(z ) \, p^{1/2}\big(z
- ( \frac{v_k}{q_k} - \frac{u_k}{q_k})\big)\,dz.
$$
and so
$$
 1 - a_k = g \big(\frac{v_k}{q_k} - \frac{u_k}{q_k} \big),
 $$
 where the function $g$ is considered in Lemma \ref{gj}.

 Using \eqref{g} and \eqref{new}, there exists $k_0$ such that,
 for any $k \ge k_0$,
 \beq
 1 - a_k \le \frac{c_{\alpha}}{2} \,
  \big|\frac{v_k}{q_k} - \frac{u_k}{q_k}
 \big|^2 \le  1/2.
 \eeq
It follows that
$$
\prod_{k \ge k_0} a_k \ge e^{- c_{\alpha} \log 2 \, (  \sum_{k \ge
k_0} \frac{|u_k - v_k|^2}{q_k^2} )\,  } >0
$$
 and so $\prod_{k \ge 1} a_k >0$.
 The second assertion follows from the first one,
 applying   \cite[Proposition
 2.19]{DZ}.
 \epf

\bre \label{zinn} {\em The result agrees with \cite[Corollary
8.1]{Zi}, which shows that
 the law of $\xi + u$, $u \in l^2$, is {\it
  absolutely continuous} with respect to the
 one of $\xi $ if and only  if
$$
 \sum_{k \ge 1} \frac{u_k^2}{q_k^2 } < \infty.
$$
 We point out that in \cite{Zi}, there are no conditions to assure
 the equivalence of $\alpha$-stable measures.
}\ere

\section{The  linear stochastic PDE }

 We start from  the linear equation
 \beq \label{ou}
dX_t  = A X_t dt  +
    dZ_t, \;\;\; x \in H.
\eeq The process   $Z$ is  a  cylindrical $\alpha$-stable process,
 $\alpha \in (0,2)$, given by
$$
 Z_t = \sum_{n \ge 1} \beta_n Z_t^n e_n,\;\;\;\; t \ge 0,
$$
 where  $(e_n)$ is the    fixed reference
 orthonormal basis in $H$,
  $(\beta_n)$ is a given sequence of {\it positive} numbers and
$(Z^n_t)$ are independent one dimensional $\alpha$-stable
 processes   defined on the same
stochastic basis $(\Omega, {\cal F}, ({\cal F}_t), \P)$,
 satisfying the usual assumptions. We have, for any $n \in \N$, $t \ge 0$,
$$
\E [ e^{i Z^n_t h}] = e^{- t|h|^{\alpha}},\;\;\; h \in \R.
$$

\begin{remark} \label{levy} {\em
 Identifying, through the basis $(e_n)$,
  the Hilbert space $H$
 with  $l^2$ and using
  Proposition \ref{ll},
  one gets that our cylindrical L\'evy process $Z$
  is a L\'evy process with values
  in the space $l^2_{\rho}$, see (\ref{rho}),
   where $(\rho_n)$ is a sequence of positive numbers
    such that
 $\sum_{n \ge 1} \beta_n^\alpha \rho_n^\alpha < \infty$.
}
\end{remark}

 We make the following assumptions.
\begin{hypothesis} \label {basic}

\hh (i)  $A : D(A) \subset H \to H$
 is a self-adjoint operator such that the fixed  basis $(e_n)$ of $H$
 verifies:
  $(e_n ) \subset D(A)$,  $A e_n  =  - \gamma_n e_n$
  with $\gamma_n >0$, for any  $n \ge 1$, and $\gamma_n \to +
  \infty $.

\hh (ii)  $\;\;\; \displaystyle{ \sum_{n \ge 1}
\frac{\beta_n^{\alpha}}
 {\gamma_n} < \infty }$ \ \ \ \
  $(\text{recall that} \;\; \beta_n > 0,\;\; \text{for any} \;\;
  n \ge 1 \;   \text{and} \; \alpha \in (0,2) ).$ \qed
\end{hypothesis}
\noindent Clearly, under (i),
 $
 D(A) = \{  x= (x_n) \in H \, :\, \sum_{n \ge 1} x_n^2
 \gamma_n^2 < +\infty
 \}.
$ In addition  $A$ generates a compact $C_0$-semigroup
  $(e^{tA})$ on $H$ such that
$$
e^{tA} e_k = e^{- \gamma_k t} e_k,\;\;
\;\; k \in \N, \;\;\; t \ge
0.
$$
 According to Hypothesis \ref{basic},
 we may consider our equation as an infinite sequence of independent
  one dimensional stochastic equations, i.e.,
 \beq \label{eq}
 dX^n_t = -
 \gamma_n X^n_t dt + \beta_n d Z^n_t,\;\;\;\; X^n _0 = x_n, \;\;
 n \in \N,
\eeq
 with   $x = (x_n) \in l^2 = H$.  The  {\it solution} is
 a stochastic process $X = (X_t^x) $ which takes  values in
   $\R^{\infty}$ with  components
 \beq \label{xt}
 X^n_t  = e^{- \gamma_n t }x_n +
  \int_0^t e^{- \, \gamma_n (t- s)
 } \beta_n dZ_s^n,\;\;\; n \in \N,\;\; t \ge 0
 \eeq

 \bpr \label{bo1}
 Assume (i) in Hypothesis \ref{basic}.
  Then, for any $x \in H$, the process $X = (X_t^x)$
   given in \eqref{xt}
   takes values in $H$ if and only if condition
  (ii) holds. Under (ii)  it can  be written as
  \begin{align} \label{mildo}  &  X_t^x = \sum_{n \ge 1} X_t^n e_n  =
e^{tA} x + Z_A(t),\;\;\; \mbox{where} \;\; \\
\nonumber &  Z_A(t) = \int_0^t e^{(t-s)A} dZ_s = \sum_{n \ge 1}
\Big( \int_0^t e^{- \, \gamma_n (t- s)
 } \beta_n dZ_s^n \Big) e_n.
  \end{align}
   The process
   $(X_t^x)$ is ${\cal F}_t$-adapted, $x \in H$. Moreover
     $X$ is  Markovian.
 \epr \noindent {\it Proof.}
  Let us consider the stochastic convolution
 \beq \label{yy}
Y_t^n=  Z_A^n(t) = \int_0^t e^{- \, \gamma_n (t- s)
 } \beta_n dZ_s^n,\;\;\; n \in \N, \;\; t \ge 0.
 \eeq
 A direct calculation shows that, for any $h \in \R,$
\beq \label{ritorna}
  \bal  \E  [e^{i h Y^n_t} ] = \exp \Big[
   -  \beta_n^{\alpha} \,
  |h|^{\alpha} \int_0^t e^{- \alpha \,
 \gamma_n \, s } ds \Big] = \exp \Big[
  - |h|^{\alpha} \, c_n^{\alpha}(t) \Big],
 \\  \mbox{where} \;\;  c_ n (t)
  =   \beta_n \, \Big(   \frac{1 - e^{- \alpha \,
  \gamma_n t}}{\alpha \, \gamma_n}
  \Big)^{1/\alpha}.
\eal
 \eeq
 It follows that
 \beq \label {serve}
  \E  [e^{i h \, Y^n_t} ] = \E  [e^{i h \,  c_n(t) \,
   L_n } ] ,  \;\;\; h \in \R,
\eeq
     where $(L_n)$
  are  independent
  $\alpha$-stable random variables having the same
   law $\mu_{\alpha}$
   (see \eqref{al}).
  Now the first assertion follows directly from Proposition
   \ref{ll}.

\vv  The property  that $(X^x_t)$ is ${\cal F}_t$-adapted is
  equivalent to
 the fact that each
  real process $\langle X_t^x, e_k \rangle$ is ${\cal F}_t$-adapted,
   for any $k \ge  1$, and this  clearly holds.

  The Markov property follows easily from the identity
 $$
 Z_A(t+h ) - e^{hA} Z_A (t) = \int_t^{t+h} e^{(t+h
 -s)A}dZ_s,\;\;\;\;\;
  t, \, h \ge 0. \qed
$$
\begin{example}\label{E1} {\em Consider
  the  following linear stochastic  heat
  equation  on $D= [0, \pi]^d$ with
  Dirichlet boundary conditions (see also
 \eqref{heat})
 \beq \label{heat1}
  \left\{
\bal dX(t, \xi)  &  = \triangle X(t, \xi) \,dt  + \,dZ(t, \xi), \;\;\; t>0,\\
      X(0, \xi) & = x(\xi), \;\;\; \xi \in D,
\\ X(t, \xi) &=0, \;\; t>0, \;\;\; \xi \in \partial D,
\eal \right. \eeq
 where $Z$ is a  cylindrical $\alpha$-stable process with
respect to the basis of  eigenfunctions of the Laplacian $\Delta$
in $H = L^2(D)$ (with  Dirichlet boundary conditions). The
eigenfunctions are
 $$
 e_j (\xi_1, \ldots, \xi_d) = (\sqrt{2/\pi})^d \sin (n_1 \xi_1)
 \cdots \sin(n_d \xi_d), \;\;\; \xi = (\xi_1, \ldots, \xi_d) \in
 \R^d,  $$
 $j=(n_1, \ldots, n_d) \in
 \N^d$.
 The corresponding  eigenvalues
 are $-\gamma_j$, where $\gamma_j =  (n_1^2 + \ldots + n_d^2 )$.
 The operator $A = \triangle $ with
  $D(A) = H^2(D) \cap H^1_0(D)$ verifies
   condition (i) in Hypothesis \ref{basic}. Moreover
    (see \cite[Section 4.4.3]{tri}) we have
 $$
D((-A)^{p/2})= \left\{ \bal H^{\alpha}(D) \cap H^1_0(D) &
 \;\;\; \text{if}
 \;\; 1< p \le 2,
\\
H^{\alpha}_0(D)  & \;\;\;\text{if}
 \;\; 1/2 < p \le 1,
 \\
H^{\alpha}(D)  & \;\;\; \text{if}
 \;\; 0< p \le 1/2.
 \eal \right.
$$
If we  identify $H$ with $l^2$  then
 $D((-A)^{p/2})$ can be identified with the weighted
 space $l^2_{\rho}$ (see \eqref{rho}) where $\rho= (\rho_j)$
 and $\rho_j = \gamma_j^{p/2}$.
 The corresponding dual spaces can be identified with
$l^2_{1/\rho}$ or with  Sobolev spaces of distributions $H^{- p}
(D)$.

  By considering sequences
  $({\beta}_j)$ of the form $({\beta}_j)= ({\gamma}_j^\delta)$ one can
 easily indicate  Sobolev spaces of distributions in which the
  cylindrical L\'evy process $Z$ might
 evolve and, at the same time,
  the Ornstein-Uhlenbeck process $X$ has
trajectories in $L^2 (D).$

For instance, assume  that $Z$ is a standard cylindrical
$\alpha$-stable process, that is ${\beta}_j = 1$,
 for any $j = (n_1, \ldots, n_d)\in \N^d$. If $p>0$
 and $\rho_j = \gamma_j^{p/2}$, then
 $l^2_{1/\rho}$ can be identified with $H^{-p}(D)$. But
 (see Proposition \ref{ll})
 $$
 \sum_{j \in \N^d} (Z_t^{j})^2 \, \gamma_j^{- p} < \infty,\;\;\;
 t>0, \;\;\;
  \P-a.s.,
$$
if and only if $\sum_{j \in \N^d} \gamma_j^{-\alpha p/ 2} <
\infty$ and if and only if $\alpha p >d$. Consequently,
 $Z_t \in H^{-p} (D)$, $t>0,$ if and only if $p > \frac{d}{\alpha}$.
} \qed
\end{example}

\subsection{Time  regularity  of trajectories}

 If the cylindrical L\'evy process $Z$ in \eqref{ou} takes values
in the Hilbert space $H$  then, by the Kotelenez regularity result
 (see \cite[Theorem 9.20]{PeZ})  trajectories of the process $X$
 which solves \eqref{ou} are c\`adl\`ag with values in $H$. However
  $Z_t\in H,$ for any $\,t
>0$, if and only if \beq
 \sum_{k \ge 1} \beta_k^{\alpha} < \infty,
 \eeq
and this is a very restrictive assumption. We conjecture that the
c\`adl\`ag property holds under much weaker conditions but, at the
moment, we are able to establish a  weaker time regularity of the
solutions.

\bth \label{kw}  Assume  Hypothesis \ref{basic}. Then the
Ornstein-Uhlenbeck process $X = (X_t^x) $
 satisfies:

\vv (i) \ for any $x \in H$, $X$  is stochastically continuous;

\vv (ii) \ for any $x \in H$, $T>0,$ $X$  has trajectories in
 $L^p(0,T; H)$, for any $0 < p< \alpha $, $\P$-a.s..
    \eth
     \bpf
 Let $0< p < \alpha$. We set $Y_t = Z_A(t)$, $t \ge 0$, and first
 show that
  \begin{equation} \label{conto}
 \E |Y_t |_{}^p  \le  \tilde c_p  \Big (  \sum_{n \ge
1} |\beta_n|^{\alpha} \frac{(1- e^{- \alpha \gamma_n t })}{\alpha
 \gamma_n }
\Big)^{p/\alpha},
 \end{equation}
  where the constant $\tilde c_p$ depends only on $p$.
 Recall that $(X_t^x)$ and
$(Y_t)$ are defined on  the same stochastic
  basis $(\Omega, {\cal F}, ({\cal F}_t)_{t \ge 0}, \P)$.
Consider a new  probability space $(\Omega', {\cal F}', \P')$
 where a Rademacher sequence $(r_n)$ is defined (i.e., $r_n
 : \Omega' \to \{1, -1\}$
  are independent and identically distributed with
   $\P' (r_n = 1)=$ $ \P' (r_n = -1) = 1/2$).

The following Khintchine inequality holds, for arbitrary real
numbers $c_1, \ldots, c_n$, for any $p>0$,
$$
 \Big( \sum_{n \ge 1} c_n^2  \Big)^{1/2} \le c_p \Big (\E' \Big|
\sum_{n \ge 1} r_n c_n  \Big|^p \Big)^{1/p},
$$
 where the constant $c_p $ depends only on $p$ (for $p=1$,
 we have $c_1 = \sqrt{2}$) and $\E'$ indicates the expectation
  with respect to $\P'$.

We fix $\omega \in \Omega$, $t \ge 0$,  and write
$$
 \Big( \sum_{n \ge 1} |Y_t^n(\omega)|^2
  \Big)^{1/2} \le c_p \Big ( \E' \Big|
\sum_{n \ge 1} r_n Y_t^n(\omega)  \Big|^p \Big)^{1/p}.
$$
Integrating with respect to $\omega$ and using the Fubini theorem
on the product space $\Omega \times \Omega'$, we find

 \begin{equation} \label{f9}
 \E |Y_t |_{}^p  \le c_p^p \,
 \E \Big [\E' \Big|
\sum_{n \ge 1} r_n Y_t^n  \Big|^p \Big]  = c_p^p
 \, \E' \Big [\E
\Big| \sum_{n \ge 1} r_n Y_t^n  \Big|^p \Big]
\end{equation}
$$
= c_p^p  \,  \E' \Big [\E \Big| \sum_{n \ge 1} r_n \int_0^t e^{-
\, \gamma_n (t- s)
 } \beta_n dZ_s^n  \Big|^p \Big].
$$
Since, for any $t \ge 0,$ $\lambda \in \R$ (using also that
$|r_n|=1$, $n \ge 1$),
 $$
 \E [ e^{i  \lambda \sum_{n \ge 1} r_n Y_t^n }] =
   e^{- |\lambda|^{\alpha} \sum_{n \ge 1}
   |\beta_n|^{\alpha} \int_0^t e^{- \alpha \, \gamma_n (t- s)
 } ds  },
$$
  we get easily
   assertion \eqref{conto}.

\vskip 1 mm \noindent (i)     It is enough to show that, for any
 $\epsilon>0$,
\begin{equation} \label{conta}
 \lim_{h \to 0^+} \, \sup_{t \ge 0}
 \P (|Y_{t+h} - Y_t| > \epsilon ) =0.
\end{equation}
Note that, for any $t \ge 0,$ $h \ge 0$,
$$
 Y_{t+h} - Y_{t} = \int_t^{t+h}  e^{(t+h -s)A}  dZ_s
 + e^{hA} \int_0^t e^{(t-s)A}  dZ_s - \int_0^t e^{(t-s)A}  dZ_s
$$
$$
= e^{hA} Y_t - Y_t + \int_t^{t+h}  e^{(t+h -s)A}  dZ_s.
$$
Let us choose $p\in (0, \alpha)$. We have
$$
 \P (|Y_{t+h} - Y_t| > \epsilon ) \le \P (\big|e^{hA}Y_{t} - Y_t \big| >
 \frac{\epsilon}{2} ) + \P ( \big | \int_t^{t+h}
  e^{(t+h -s)A} dZ_s \big| >
 \frac{\epsilon}{2} )
$$
$$
\le 2^{p} \frac{\, \E |e^{hA}Y_{t} - Y_t  |^p}{\epsilon^p}
 +
 2^{p}\frac{ \, \E | \int_0^{h}
  e^{sA} dZ_s \big|^p}{\epsilon^p} = I_1(t,h) + I_2(h).
$$
But (see \eqref{conto})
 $$
 \E |Y_t |_{}^p  \le  c_p  \Big (  \sum_{n \ge
1} |\beta_n|^{\alpha} \frac{(1- e^{- \alpha \gamma_n t })}{\alpha
 \gamma_n }
\Big)^{p/\alpha}
$$
and so
$$
 [I_2(h)]^{\alpha/p} \to 0,\;\;\; \mbox{as} \;\; h \to 0^+.
$$
Concerning $I_1$, we find, using again the Khintchine inequality,
$$
|e^{hA} Y_t - Y_t|  =  \Big( \sum_{n \ge 1} |(e^{- \gamma_n h }
-1)Y_t^n|^2 \Big)^{1/2} \le c_p \Big ( \E' \Big| \sum_{n \ge 1}
r_n (e^{-\gamma_n h} -1)Y_t^n \Big|^p \Big)^{1/p}
$$
and, reasoning as in \eqref{f9} with $\beta_n$ replaced by
 $(1 - e^{- \gamma_n h})\beta_n$,
$$
\E |e^{hA} Y_t - Y_t|^p  \le  c_p^p \,   \E' \E \Big| \sum_{n \ge
1} r_n (e^{-\gamma_n h} -1)Y_t^n \Big|^p $$ $$\le C_p \Big (
\sum_{n \ge 1} |(1 - e^{- \gamma_n h})\beta_n|^{\alpha} \frac{(1-
e^{- \alpha \gamma_n t })}{\alpha
 \gamma_n }
\Big)^{p/\alpha} \le \frac{C_p}{\alpha^{p/\alpha}} \Big (  \sum_{n
\ge 1} \frac{|(1 - e^{- \gamma_n h})\beta_n|^{\alpha}} {
 \gamma_n }
\Big)^{p/\alpha},
$$
 $t \ge 0.$
Since
$$
\lim_{h \to 0^+} \Big (  \sum_{n \ge 1} \frac{|(1 - e^{- \gamma_n
h})\beta_n|^{\alpha}} {
 \gamma_n }
\Big)^{p/\alpha} =0,
$$
we get
$$
\lim_{h \to 0^+} \, \sup_{t \ge 0} I_1 (t,h) =0
$$
and so assertion \eqref{conta} is proved.

\hh {(ii) }   It is enough to show that
 \beq
\E  \int_0^T \Big( \sum_{n \ge 1} |Y_t^n|^2 \Big)^{p/2} dt <
\infty,
  \eeq
where $Y_t = Z_A(t)$, $t \ge 0$. Using \eqref{conto} we get
$$
\int_0^T \E |Y_t |_{}^p dt \le \tilde c_p \int_0^T \Big (  \sum_{n
\ge 1} |\beta_n|^{\alpha} \frac{(1- e^{- \alpha \gamma_n t
})}{\alpha
 \gamma_n }
\Big)^{p/\alpha} dt \le C_{p,\alpha} \, T \Big ( \sum_{n \ge 1}
\frac{ |\beta_n|^{\alpha}} {\gamma_n } \Big)^{p/ \alpha} <
+\infty.
$$
The proof is complete.
 \epf

\begin{remark} {\em  In the limiting
Gaussian case of $\alpha =2$,  the previous proof allows to
  get
the well known result that
 trajectories of $X$ are in $L^2(0,T; H)$, for any $T>0.$
}
\end{remark}


  Using  that $X = (X_t^x)$, $x \in H$, is stochastically
 continuous and ${\cal F}_t$-adapted (see Theorem \ref{kw} and
 Proposition \ref{bo}) we can apply \cite[Proposition 3.6]{DZ} and
 obtain
\begin{corollary}\label{pred} For any $x \in H$, the process
   $(X_t^x)$   has a  predictable version.
\end{corollary}
\noindent For $p \in (0,1)$, $L^p(0,T; H)$ is a
 linear complete and separable metric
   space with respect to the distance
 $
 d_p (f,g) = \int_0^T |f(t) - g(t)|^p dt , $ $
  f,g \in L^p(0,T; H).
$ From Theorem \ref{kw} it is  straightforward to obtain
 \bco \label{kw1} Assume  Hypothesis \ref{basic}. Then, for any $T>0$,
  $x \in H$,  $\P$-a.s.,
 the Ornstein-Uhlenbeck process $X = (X_t^x)_{t \in [0,T]} $ is a
 random variable with
 values
 in
 $L^p(0,T; H)$,   for any $0 < p< \alpha $.
\eco

\subsection{Support}

 We start with a preliminary one dimensional result.

\bpr \label{dopo}
 Let $L= (L_t)$ be a one dimensional $\alpha$-stable process, $\alpha \in
 (0,2)$. Let  $\gamma \in \R $  and set
 \begin{equation} \label{kt}
K(t) = \int_{0}^{t} e^{\gamma(t-s)} dL_s,\,\,\, \;\;\; t\geq 0.
 \end{equation}
Then, for any $p >0$, $T>0$, the random variable
  $(K, K_T)$ has full support in $L^p(0,T) \times \R$.
\epr The proposition is a direct corollary of the following
general lemma.

\begin{lemma} \label{le}
Let $L=(L_t)$ be a real valued
 L{\'e}vy process with intensity
 measure $\nu$ (see (\ref{psi})).
 Suppose that there exists $R>0$
  such  that $\nu$ restricted to $(-R, R)$ has an
 absolutely continuous
 part with a strictly positive density.
 Then, for any $p >0$, $T>0$, the random variable
  $(L, L_T)$ has full support in $L^p(0,T) \times \R$.
\end{lemma}
\bpf We set $\nu = 1_{(-R,R)} \nu +  1_{ \R \setminus (-R, R) }
\nu = \mu_0 + \mu_1$, where $ \mu_0 = 1_{(-R,R)} \nu $ denotes the
Borel measure such that
 $ \mu_0 (A) = \nu (A \cap (-R, R))$, for any Borel set
 $A \subset \R$.  Let $g$ be the density of the absolutely
  continuous  part
  $[\mu_0]_{ac}$ of $\mu_0$  and let $[\mu_0]_{s}$ be the singular
 part of $\mu_0$.
Then
$$
\mu_0 = [\mu_0]_{ac} + [\mu_0]_{s} =
 (g \wedge 1)1_{(-R,R)} {\cal L}_1 \, + \,
 (g - [g \wedge 1] )1_{(-R,R)} {\cal L}_1 \, +  \, [\mu_0]_{s}
$$
(${\cal L}_1$ denotes the one dimensional Lebesgue measure). We
write  $\nu = \nu_0 + \nu_1$, where
$$
 \nu_0 = \mu_1  \, + \,
 (g - [g \wedge 1] )1_{(-R,R)} {\cal L}_1 \, +  \, [\mu_0]_{s}
 \;\;\; \mbox{and} \;\;
\nu_1 = (g \wedge 1)1_{(-R,R)} {\cal L}_1.
$$
Note that  $\nu_0$ and  $\nu_1$ are both positive measures and
moreover ${\nu_1}$ is a finite measure with a  positive density $g
\wedge 1$ on $(-R, R)$.

 Let us introduce two independent L\'evy processes $L^0 = (L^0_t) $
 and $L^1 = (L^1_t)$.  The exponent $\psi$
  of $L^0$ (see \eqref{psi})
   is the same
 of  $L $ but with the jump intensity measure $\nu$ replaced by
 $\nu_0$. The exponent of $L^1 $ is given in \eqref{psi}
  with $Q =0$, $a=0$ and
  $\nu$ replaced by $\nu_1$, i.e.,
  $L^1$ is a compound L\'evy process with
   intensity measure $\nu_1$.

By using the characteristic function and independence, we obtain
that the process
$$
 \tilde L = L^1 + L^0,\;\;\; \text{i.e.,} \;\; \tilde
 L_t = L^1_t + L^0_t,\;\; t
\ge 0,
$$
 has the same law of the initial L\'evy process $L$.
  It follows that  the law of $(L, L_T)$ is the convolution of the laws of
    $(L^0, L_T^0)$ and  $(L^1, L_T^1)$.   Our assertion will follow
 from the fact that
 $(L^1, L_T^1)$ has full support in $L^p(0,T) \times \R$.

 Taking into account that pice-wise constant
 functions taking value 0 at $t=0$
  are dense in $L^p(0,T)
 $, for any $p>0$, we only have to prove that for a
 fixed pice-wise constant function
  $\phi : [0,T] \to \R$, with $\phi (0)=0$,
    for a fixed $a \in \R$ and $\epsilon >0$,
 \beq \label{t}
 \P \Big ( \int_0^T |L^1_t - \phi(t)|^p dt + |L^1_T -a| < \epsilon
  \Big ) >0.
\eeq We may  assume that $\phi(T)=a$ and that  $\phi$ takes real
values $0,x_1, \ldots, x_{k-1}$, $x_k=a$,
 respectively on intervals $[0, t_1[$, $\ldots,$ $[t_k, T[$, with
  $0 < t_1< \ldots $ $t_k< T$. Define
$$
S = \sup \{ |x_i|,\;\;\; i=1, \ldots, k \}.
$$
Let $0< \tau_1 <$ $\ldots < \tau_k $ be the first $k$ consecutive
moments of jumps for the process $L^1$ and denote by $Y_1$,
$\ldots, Y_k$
 the random variables $L^1_{\tau_1}$, $\ldots, L^1_{\tau_k}$;
 set $Y_0=0$ and $\tau_0 =0$.

Note that $\tau_j - \tau_{j-1}$, $j =1, \ldots, k ,$ and
 $Y_j - Y_{j-1}$, $j =1, \ldots, k ,$
  are independent
  random variables. Moreover,
 $\tau_j - \tau_{j-1}$  have the same exponential distribution and
 $Y_j - Y_{j-1}$ have the   positive density $g \wedge 1 \, \big(
  \int_{(-R,R)} \,
   g \wedge 1 \big)^{-1}$
  on $(-R, R)$.

\vskip 1mm  For arbitrary $i,j \in \{ 0, \ldots , k\}$, $\delta
 >0$, $M > S-R $
  the independent events
$$
\{ |\tau_i - t_j | \le  \delta \}, \;\;\; \{ |Y_i - x_j| \le M \}
$$
have all positive probabilities. Using this fact and the property of
independence, we get easily \eqref{t}. \epf

\vv {\bf Proof of Proposition \ref{dopo}.} We consider  $\gamma
\not =0$  (the case $\gamma =0$ follows from Lemma \ref{le}).
Using \cite[Theorem 3.1]{RWo}, we know that there exists an
$\alpha$-stable process
 $Z = (Z_t)$ such that
$$
\int_{0}^{t} e^{- \gamma s} dL_s = Z(h(t)), \;\;\; \text{where}
 \;\; h(t)=\frac{1 - e^{- \alpha \gamma t}  }{ \alpha \gamma},
 \;\; \;\; t \ge 0.
$$
 Consequently,
 $ K(t) =  e^{\gamma t} Z(h(t))$ (see \eqref{kt}).
   Using Lemma \ref{le} and the
 fact that $h \in C^{\infty}([0, +\infty[)$
  with $h'(t) \not = 0 $, $t \ge 0$, we get easily the assertion.
\qed

\bth \label{bo}
  Assume  Hypothesis \ref{basic} and fix $T>0$, $x \in H$
   and $p \in
   (0, \alpha)$.
     Consider  the Ornstein-Uhlenbeck
    process $X = (X_t^x)_{t \in [0,T]} $, solving \eqref{ou}.
     The support of the random variable
     $(X, X_T^x): \Omega \to L^p(0,T; H) \times H$
      is $L^p(0,T; H) \times H$.
 \eth
\bpf Let $X^x_t = (X_t^n)$, $t \ge 0$.
 It is enough to prove that, for any $\epsilon
>0$, and for any $(\phi, a) \in L^p(0,T; H) \times H$,
 one has
$$
 \P \Big ( \int_0^T  \Big ( \sum_{n \ge 1 }
 |X_t^n - \phi_n(t) |^2   \Big)^{p/2} dt  < \epsilon, \;
\sum_{n \ge 1} | X^n_T - a_n |^2 < \epsilon \Big) >0.
$$
By using a standard density argument, we may assume that
  $(\phi, a)$ is of the form
 $$
\phi(t) = \sum_{k=1}^N \phi_k(t) e_k,\;\;\; a = \sum_{k=1}^N a_k
e_k,
$$
for some $N \in \N.$ We write, using that $p/2<1$,
 $$
\P \Big ( \int_0^T  \Big ( \sum_{n = 1 }^N
 |X_t^n - \phi_n(t) |^2   \Big)^{p/2} dt  < \epsilon, \;
\sum_{n = 1}^N | X^n_T - a_n |^2 < \epsilon \Big)
$$
$$
\ge \P \Big ( \int_0^T  \sum_{n=1}^N
 |X_t^n - \phi_n(t) |^p   dt  < \epsilon, \;
\sum_{n = 1}^N | X^n_T - a_n |^2 < \epsilon \Big)
$$
$$
\ge  \P \Big ( \int_0^T
 |X_t^1 - \phi_1(t) |^p   dt  < \epsilon/N, \;
 | X^1_T - a_1 |^2 < \epsilon / N \Big) \cdots
 $$ $$ \cdots  \P \Big ( \int_0^T
 |X_t^N - \phi_N(t) |^p   dt  < \epsilon/N, \;
 | X^N_T - a_N |^2 < \epsilon/N  \Big),
$$
 using independence. By Proposition \ref{dopo} we know that
 the previous product of probabilities is positive. The proof is
 complete.
\epf

\bco \label{irro}
 Under  Hypothesis \ref{basic}, for any $x \in
H$, the OU
 process $(X_t^x)$
  is irreducible, i.e., for any open ball $B \subset H$, $t>0,$
 we have $ \P (X^x_t \in B)
>0.
$ \eco

\subsection{Equivalence of transition probabilities}

Here we will  assume Hypothesis \ref{basic} together with
\begin{hypothesis} \label{basic1} For any $t>0,$
\begin{equation} \label{e4}
  \;\;  \sup_{n \ge 1} \frac{ e^{-  \gamma_n t }\, \gamma_n
^{1/ \alpha}} { \beta_n^{}} = C_t
 < \infty.
 \end{equation}
\end{hypothesis}

\bth \label{ass} Assume  Hypotheses \ref{basic} and \ref{basic1}.
   Then the laws   $\mu^x_t$
  and $\mu^y_t$  of $X_t^x$ and $X_t^y$, respectively,
   are equivalent, for any
  $t>0,$  $x, y \in H$, $\alpha \in (0,2)$.
Moreover  the density $\frac{d \mu_t^x}{d \mu_t^y}$ of
 $\mu_t^x$ with respect to $\mu_t^y$ is given by
 $$
\frac{d \mu_t^x}{d \mu_t^y} = \lim _{N \to \infty} \prod_{k=1}^N
  \frac{p_{\alpha} \Big(\frac{ z_k - e^{- \gamma_k t }x_k}{c_k (t)}
 \Big)}{p_{\alpha} \Big(\frac{ z_k - e^{- \gamma_k t }y_k}{c_k(t)}
 \Big)} \;\;\; in \;\; L^1 (\mu_t^y), \;\;  \mbox{where}
 \;\;  c_ k (t) =
   \beta_k \, \Big(   \frac{1 - e^{- \alpha \,
  \gamma_k t}}{\alpha \, \gamma_k}
  \Big)^{1/\alpha},
 $$
 and $p_{\alpha}$ is the density  of
  the one dimensional  $\alpha$-stable measure
   considered in (\ref{al}).

 If \eqref{e4} does not hold then for some $x \in H$, $\mu^x_t$ is
 not  absolutely continuous with respect to $\mu^0_t$ (the law of
 $X_t^0 = Z_A(t), see \ (\ref{mildo}))$.
 \eth

\begin{remark} \label{df}
 {\em If we assume  Hypothesis \ref{basic},
 then Hypothesis \ref{basic1}
   is  sharp
  in the limiting Gaussian  case of $\alpha=2$. Indeed, under
   Hypothesis \ref{basic} and $\alpha=2$,
    Hypothesis \ref{basic1}  is equivalent to each of the following
     facts:

\vv (i) \ the laws of  $X_t^x$ and $X_t^y$ are equivalent,
 for any
 $t>0$,  $x, y \in H$;

\vv (ii) \      the
  Gaussian Ornstein-Uhlenbeck semigroup $(R_t)$
 associated to $(X_t^x)$
  is strong Feller    (see  \cite[Section 9.4.1]{DZ}).

\vv In addition, under   Hypothesis \ref{basic} and $\alpha=2$,
the following regularizing property:
 $$
 R_t f \in C^{\infty}_b (H), \;\;  t>0, \;\; f \in B_b (H),
$$
 holds if and only  if $( e^{-  \gamma_n t }\,
 \sqrt{ \frac{\gamma_n ^{}} { \beta_n^{2}} } )$
  is a bounded sequence.
}
\end{remark}

 \bpf [\textbf{Proof of Theorem \ref{ass}.}] Fix $x = (x_n)$ and $y=(y_n)$.
  Let $Y_t = Z_A(t)$ and  $p = p_{\alpha}$.
 Consider formulas \eqref{xt}, \eqref{yy} and \eqref{serve}.
  The density of the random variable $Y_t^k$
    is clearly
 $\frac{1}{c_k (t)} p \Big(\frac{ z_k} {c_k(t)} \Big)$  so that the
 density of $X_t^k$ is
$
 \frac{1}{c_k (t)} p \Big(\frac{ z_k - e^{- \gamma_k t }x_k} {c_k(t)}
 \Big).
$
 The measures  $\mu_t^x$ and $\mu_t^y$ can be seen as Borel
  product measures
 in $\R^{\infty}$, i.e.,
$$
 \bal
& \mu_t^x = \prod_{k \ge 1} \mu_t^{x_k},\;\;\;\;\;
   \mu_t^y = \prod_{k \ge 1} \mu_t^{y_k}, \;\;\; \mbox{where}
  \;\;  \mu_t^{x_k}, \; \mu_t^{y_k} \;\;  \mbox{have densities, respectively, }
\\
&  \frac{1}{c_k (t)} p \Big(\frac{ z_k - e^{- \gamma_k t }x_k}
{c_k(t)}
 \Big)  \;\; \mbox{and} \;\;
  \frac{1}{c_k (t)} p \Big(\frac{ z_k - e^{- \gamma_k t }y_k} {c_k(t)}
 \Big).
 \eal
$$
 To get the assertion we will apply  Theorem \ref{gh}.
 To this purpose, one
  checks that
$$
\sum_{k \ge 1} \frac{e^{- 2\gamma_k t} { |x_k - y_k|^2}
 }{c_k (t)^2} <
 \infty.
$$
This follows easily from \eqref{e4}.

 If
 \eqref{e4} does not hold, for some $t>0$,
   then it is easy to see that there exists
  $\hat x = (\hat x_n) \in H$ such that
$$
 \sum_{k \ge 1} \frac{e^{- 2\gamma_k t} { \, \hat x_k^2}
 }{c_k (t)^2} =
 + \infty.
$$
According to Remark \ref{zinn}, this condition means that
$\mu_t^{\hat x}$, the law  of $X^{\hat x}_t$, is not absolutely
continuous with respect to $\mu_t^0$. \epf

\subsection{Smoothing effect}

We now consider the transition Markov semigroup $(R_t)$ associated
to $(X_t^x)$, i.e.
 $R_t : B_b (H) \to B_b (H)$,
$$
 R_t f (x) = \E [f(X_t^x)],\;\;\; x\in H,\;\; f \in B_b(H),\;\;\;
 t \ge 0.
$$
 The next result shows    that $(R_t)$
  has a smoothing effect and
   that gradient estimates hold for it.
   In particular, this will imply the strong
   Feller property for $(R_t)$.
    Recall that  a   Markov semigroup
 $(P_t) $ acting on
  $B_b (H)$ is said to be  {\it strong Feller,} if
  $ P_t f \in C_b (H),$  for any $ t>0$  and
  $f \in B_b(H).$

\bth \label{grad} Assume Hypotheses \ref{basic} and \ref{basic1}.
 Then, for any $t>0,$ the transition semigroup $(R_t)$ maps Borel and bounded
 functions
 into $C^1_b(H)-$functions.  Moreover, for any $k \in H$
 with $|k| \le 1$, $f \in B_b(H),$ $t>0,$
   we have
 \beq
 \sup_{x \in H}|\lan DR_t f(x),
 k \ran | \le 8 c_{\alpha} \,   \tilde C_t  \| f\|_{0},
 \;\;\;  \mbox{ where} \;\;\;
 \tilde C_t =  \sup_{n \ge 1} \frac{
\gamma_n ^{1/ \alpha} \alpha^{1/\alpha}} { \beta_n^{} \, (e^{\alpha \gamma_n t} -1)^{1/\alpha}}
\eeq
 ($c_{\alpha}$ is defined in (\ref{g})).
 Finally, for any $t>0$, $f \in C_b(H)$,
   $x= (x_n),$ $ h=(h_n) \in H$, we have
 \beq \label{gr}
 \lan DR_t f (x), h \ran  =
    \int_H f(e^{tA} x +  y ) \, \, \sum_{k \ge 1} \,
  \frac{p'_{\alpha}(\frac{y_k} {c_k (t)})}
  {p_{\alpha}(\frac{y_k} {c_k (t)})} \, \,
\frac{e^{- \gamma_k t  } h_k}{c_k (t)} \,
  \; \mu^0_t (dy),
\eeq where
 $\mu^0_t$ is the law of $X_t^0 = Z_A(t)$.
 \eth
\bre {\em Set $\gamma_{min} = \inf_{n \ge 1} \gamma_n >0$. We have, for any $t>0$,
\begin{align} \label{ritorno}
\tilde C_t \le \alpha^{1/\alpha} \sup_{s \ge  t\gamma_{min}}
\frac{ e^s} {
\, (e^{\alpha s} -1)^{1/\alpha}} \, \cdot \,
\sup_{n \ge 1} \frac{
\gamma_n ^{1/ \alpha} e^{-\gamma_n t}} { \beta_n^{}
\, } \le c(t) \, \cdot C_t
\end{align}
and $C_t < \infty$ by Hypothesis \ref{basic1}.
}\ere

\bpf  We fix $t>0$.

   The proof is divided into some
  steps. By the first three steps, we will  show that, for any
   $f \in C_b(H)$,
   $R_t f$ is G\^ateaux differentiable at any $x \in H$
  and moreover that equality \eqref{gr} holds.

 \vv {\it I Step.} \
   We  assume that $f \in C_b (H)$ is  cylindrical, i.e., it
 depends only  on a finite numbers of coordinates. Identifying
  $H$ with $l^2$ through the basis $(e_n)$, we  have
 \begin{equation} \label{cy}
f (x) = \tilde f(x_1, \ldots,
  x_j), \;\;\; x \in H,
  \end{equation}
 for some  $j \ge 1$, and  $\tilde f : \R^j \to \R$ continuous and bounded.
  In this
 first step we also assume that
   $\tilde f$ has bounded support in $\R^j$.

Fix  arbitrary $x, h \in H$. We want to show that
 there exists $D_h R_t f (x)$, the directional derivative of $R_t
  f$
  at $x $, along the direction $h$.
   Set $h_N = \sum_{k =1}^N h_k e_k$
   so that $h_N \to h $ in $H$.
  Since $f$ is cylindrical, for $m \ge \max(j,N)$, we get
$$
 R_t f (x) = \int_{H } f(y) \prod_{k \ge 1}  \mu_t^{x_k} (dy)
 = \int_{\R^{m}} \tilde f(z)\,
  \prod_{k=1}^m
  p_{\alpha} \Big(\frac{ z_k - e^{- \gamma_k t }x_k}{c_k (t)}
 \Big) \frac{1}{c_k (t) } \, dz_k.
$$
 Using our assumptions on $\tilde f$, it is not difficult to show that
 there exists
 \beqr
D_{h_N} R_t f (x) =  - \int_{\R^{m}} \tilde f(z)\,
  \Big ( \sum_{k=1}^N \frac{
  p_{\alpha}' \Big(\frac{ z_k - e^{- \gamma_k t }x_k}{c_k (t)}
 \Big)} {p_{\alpha}
 \Big(\frac{ z_k - e^{- \gamma_k t }x_k}{c_k (t)}
 \Big)} \frac{ e^{- \gamma_k t}  h_k  }{c_k (t)} \Big) \,
 \cdot
 \; \\  \cdot \, \prod_{k=1}^m
  p_{\alpha} \Big(\frac{ z_k - e^{- \gamma_k t }x_k}{c_k (t)}
 \Big) \frac{1}{c_k (t) } \, dz_k
\\
=  - \int_{H} f(z)\,
  \Big ( \sum_{k=1}^N \frac{
  p_{\alpha}' \Big(\frac{ z_k - e^{- \gamma_k t }x_k}{c_k (t)}
 \Big)} {p_{\alpha} \Big(\frac{ z_k - e^{- \gamma_k t }x_k}{c_k (t)}
 \Big)}  \frac{ e^{- \gamma_k t}  h_k  }{c_k (t)} \Big)
  \, \prod_{k \ge 1} \mu_t^{x_k} (dz_k),\;\;\; N \in \N.
\eeqr
 In order to pass to the limit, as $N \to \infty$, we show that
\beq \label{l2}
  g_N (t,x) = \sum_{k=1}^N \frac{
  p_{\alpha}' \Big(\frac{ z_k - e^{- \gamma_k t }x_k}{c_k (t)}
 \Big)} {p_{\alpha} \Big(\frac{ z_k - e^{- \gamma_k t }x_k}{c_k (t)}
 \Big)}  \frac{ e^{- \gamma_k t}  h_k  }{c_k (t)} \;\;
 \mbox{converges in} \;\; L^2 (\mu_t^x).
 \eeq
 Using that, for $j \not = k$,
 \beqr
\frac{ e^{- \gamma_k t }   h_k} {c_k (t)} \frac{ e^{- \gamma_j t}
 h_j } {c_j (t)} \int_{\R^2}
 \frac{
  p_{\alpha}' \Big(\frac{ z_k - e^{- \gamma_k t }x_k}{c_k (t)}
 \Big)} {p_{\alpha} \Big(\frac{ z_k - e^{- \gamma_k t }x_k}{c_k (t)}
 \Big)}
 \, \frac{
  p_{\alpha}' \Big(\frac{ z_j - e^{- \gamma_j t }x_j}{c_j (t)}
 \Big)} {p_{\alpha} \Big(\frac{ z_j - e^{- \gamma_j t }x_j}{c_j (t)}
 \Big)} \; \cdot
  \\
\, \cdot  p_{\alpha} \Big(\frac{ z_k - e^{- \gamma_k t }x_k}{c_k
(t)}
 \Big ) p_{\alpha} \Big(\frac{ z_j - e^{- \gamma_j t }x_j}{c_j
 (t)}\Big)
dz_k  dz_j =0
  \eeqr
 (since $p'_{\alpha}$ is odd) we get, for $N, p \in \N$,
 \beqr
 \int_{H}  \Big | \sum_{k=N}^{N +p} \frac{
  p_{\alpha}' \Big(\frac{ z_k - e^{- \gamma_k t }x_k}{c_k (t)}
 \Big)} {p_{\alpha} \Big(\frac{ z_k - e^{- \gamma_k t }x_k}{c_k (t)}
 \Big)}  \frac{ e^{- \gamma_k t}  h_k  }{c_k (t)} \Big|^2
 \mu_t^x (dz)
\\
\int_{\R^{p+1}}  \Big | \sum_{k=N}^{N +p} \frac{
  p_{\alpha}' \Big(\frac{ z_k - e^{- \gamma_k t }x_k}{c_k (t)}
 \Big)} {p_{\alpha} \Big(\frac{ z_k - e^{- \gamma_k t }x_k}{c_k (t)}
 \Big)}  \frac{ e^{- \gamma_k t } h_k  }{c_k (t)} \Big|^2
 \prod_{k =N}^{N+p} p_{\alpha} \Big(\frac{ z_k - e^{- \gamma_k t }x_k}{c_k (t)}
 \Big) \frac{1}{c_k (t) } \, dz_k
\\
= \int_{H}   \sum_{k=N}^{N +p} \frac{
  (p_{\alpha}')^2 \Big(\frac{ z_k - e^{- \gamma_k t }x_k}{c_k (t)}
 \Big)} {p_{\alpha}^2 \Big(\frac{ z_k - e^{- \gamma_k t }x_k}{c_k (t)}
 \Big)}  \frac{ e^{- 2\gamma_k t}  h_k^2  }{c_k ^2 (t)}
 \; \mu_t^x (dz)
\\
=  \sum_{k=N}^{N +p}   \frac{ e^{- 2\gamma_k t } h_k^2  }{c_k ^2
(t)} \int_{\R} \frac{
  p_{\alpha}'^2 (y_k) }
  {p_{\alpha} (y_k)}  \, d y_k \, \le \, 8 c_{\alpha}  \tilde C_t^2  \sum_{k=N}^{N +p}
  h_k^2,
\eeqr
 where $ 8 c_{\alpha} =\int_{\R} \frac{
  p_{\alpha}'^2 (y) }
  {p_{\alpha} (y)}  \, d y $ (see \eqref{g}). This proves
  \eqref{l2}.

Note that, for any $N \in \N$,
 $$
D_{h_N} R_t f (x) =  - \int_{H} f(z + e^{tA} x)\,
  \Big ( \sum_{k=1}^N \frac{
  p_{\alpha}' \Big(\frac{ z_k }{c_k (t)}
 \Big)} {p_{\alpha} \Big(\frac{ z_k }{c_k (t)}
 \Big)}  \frac{ e^{- \gamma_k t } h_k  }{c_k (t)} \Big)
  \,  \mu_t^{0} (dz).
$$
Up to now we have showed that
 \beq \label{fi}
\frac{ R_t f (x + s h_N ) - R_t f(x)}{s} = \frac{1}{s} \int_0^s
D_{h_N} R_t
 f(x + r h_N)dr, \;\; s \in (-1,1).
\eeq
 Using also \eqref{l2}, it is not difficult to show
 that, for any $r \in (-1,1)$, $N \in \N$,
  \beq \label{ri}
 \lim_{N \to \infty} D_{h_N} R_t
 f(x + r h_N) = - \int_{H} f(z + e^{tA} (x + r h) )\,
  \Big ( \sum_{k=1}^{\infty} \frac{
  p_{\alpha}' \Big(\frac{ z_k }{c_k (t)}
 \Big)} {p_{\alpha} \Big(\frac{ z_k }{c_k (t)}
 \Big)}  \frac{ e^{- \gamma_k t } h_k  }{c_k (t)} \Big)
  \,  \mu_t^{0} (dz).
 \eeq
 Moreover, $| D_{h_N} R_t
 f(x + r h_N)| \le \, 8 c_{\alpha} C_t |h| \| f\|_0$, for any $r \in
 (-1,1)$. Thus we can pass to the limit, as $N\to \infty$, in
 \eqref{fi} and get
 \beq  \label{gat}
\frac{ R_t f (x + s h ) - R_t f(x)}{s} = \frac{1}{s} \int_0^s
u(t,x + r h)\, dr, \;\; s \in (-1,1),
 \eeq
 where $u(t,x+ rh)$ is the right-hand side of \eqref{ri}.
 This shows that $R_t f$ is G\^ateaux differentiable at $x \in H$
  along
  the direction $h$ and moreover that \eqref{gr} holds.

\hh  {\it II Step.} \ We consider $f \in C_b (H)$
  which is only cylindrical (i.e., $f$ is  given by \eqref{cy}
   but the function
  $\tilde f $  is not assumed to have  bounded
  support in $\R^j$).

 \vskip 1mm Define $ \tilde f_n (y) = \tilde f(y)
\phi(\frac{|y|}{n})$, for any $y \in \R^j$,
  where $\phi : [0, +
\infty[ \to \R_+$ is a continuous function such that, $\phi (s)
=1$, $s \in [0,1]$, $\phi (s) =0, $ $s \ge 2$.

 \vskip 1mm We have that
 $\| \tilde  f_n \|_0 \le  \| \tilde f\|_0$, $n \in \N$,
  and moreover
$ \tilde  f_n (y) \to \tilde f(y)$, as $n \to \infty$, for any $y
\in \R^j$.

Let $f_n : H \to \R$, $f_n (x) = \tilde f_n (x_1, \ldots, x_j) $,
 for any $x \in H$, $n \in \N$.

 We find by the previous step, for any $n \in \N$ and $x \in H$,
 \beq
\frac{ R_t f_n (x + s h ) - R_t f_n(x)}{s} = \frac{1}{s} \int_0^s
D_{h} R_t
 f_n (x + r h)dr, \;\; \;\;  s \in (-1,1).
 \eeq
 Passing to the limit, as $n \to \infty$, it is easy to see
  that \eqref{gat}
  holds for $f$. This shows the G\^ateaux differentiability of $R_t f
   $ on $H$ and also the equality \eqref{gr}.

\hh  {\it III Step.} \ We consider  an arbitrary $f \in C_b (H)$.
Let us introduce  the cylindrical functions $g_n,$ $
 g_n (x) = f \Big( \sum_{k =1}^n x_k e_k \Big),
 $ $ n \in \N, \;\; x \in
 H.
$

 It is clear that $\| g_n \|_0 \le  \| f\|_0$, $n \in
\N$, and moreover $g_n (x) \to f(x)$, for any $x \in H$. Repeating
the  argument of the previous step, with $f_n$ replaced by $g_n$,
 and passing to the limit, we get that the assertion of the previous
  step holds even  for any $f \in C_b (H)$.

\hh {\it IV Step.} \ Let $f \in  C_b (H) $ and consider the
G\^ateaux derivative of $R_t f$ in $x \in H$
$$
 DR_t f (x) = \int_H f(e^{tA} x +  y ) \,  \sum_{k \ge 1} \,
  \frac{p'_{\alpha}(\frac{y_k} {c_k (t)})}
  {p_{\alpha}(\frac{y_k} {c_k (t)})} \,
\frac{e^{- \gamma_k t }}{c_k (t)} \, e_k
  \; \mu^0_t (dy).
$$
It is not difficult to show that $D R_t f : H \to H$ is
continuous. This gives that {\sl  $R_t f$ is Fr\'echet
differentiable at any $x \in H$.   }
   Moreover, we have  the required gradient estimate
$$ \| DR_t f \|_0 \le 8 c_{\alpha}   \tilde C_t  \| f
\|_0
$$
 {\it V Step.} \ To complete the proof, take $g \in B_b (H)$.
  A well known argument
   (see  \cite[Lemma 7.1.5]{DZ1})
  shows that  $R_t g $ is Lipschitz
 continuous on $H$, for any $t>0$.
  Then  the semigroup law  gives that
 $R_t g \in {
C}_b^1(H),$  for any $t>0$.
 The proof is complete.
  \epf

\bre {\em Under the assumptions of Theorem \ref{grad}, one could
 show   the following regularizing property
 $
 R_t f \in C^{\infty}_b (H), $ $  t>0,$ $ f \in B_b (H).
$
  This  generalizes the well  known smoothing
 property of the Gaussian Ornstein-Uhlenbeck
 semigroup (see  Remark \ref{df}).
} \ere

\bre{ \em  Theorem \ref{ass} can be also deduced  from Theorem
\ref{grad} and Corollary \ref{irro} if one applies the Hasminkii
theorem (see \cite[Proposition 4.1.1]{DZ1}).
 } \ere

\section{Nonlinear stochastic PDEs}

 We pass now to
 nonlinear
  SPDEs of the form
 \beq \label{ab1}
 dX_t  = A X_t dt + F(X_t)dt +  dZ_t, \;\;\;\; X_0  = x \in l^2 = H,
 \eeq
where $Z = (Z_t)$ is a cylindrical $\alpha$-stable
 L\'evy
 process.
 Throughout the section,  we will assume Hypothesis \ref{basic}
  and also
\begin{hypothesis} \label {dl}
$ \;\; F: H \to H \;\; \mbox{is Lipschitz continuous and
 bounded.}$
\end{hypothesis}

\subsection{Existence, uniqueness and Markov property}

We say that a predictable
   $H$-valued
   stochastic process $X= (X_t^x)$, depending on  $x
 \in H$,   is
 a {\it mild solution} to equation
\eqref{ab1} if, for any $t \ge 0,$ $x \in H$, it holds: \beq
\label{mil}
  \bal
 X_t^x = e^{tA} x + \int_0^t e^{(t-s)A} F(X_s^x)ds
   + Z_A (t),\;\; \P-a.s., \;\;\;
 \mbox{where} \;\;
Z_A(t) = \int_0^t e^{(t-s)A}  d Z_s,
 \eal
 \eeq
see \eqref{mildo}. In  formula \eqref{mil}
 we are considering  a predictable version of the process
  $(Z_A (t))$ according to Corollary \ref{pred}.

Note that, since $F$ is bounded,  the deterministic integral in
\eqref{mil} is a well defined continuous process. Moreover,
  as far as the regularity of trajectories is concerned,
 the mild
 solution will have
 the same regularity  as
  $(Z_A(t))$.
  In particular, according to Theorem
  \ref{kw}, any mild solution $X$ will be stochastically
   continuous.

 \vskip 1mm  To show existence and uniqueness we need the
 following deterministic result which is not standard in
   the case $p \in (0,1)$.

   Recall that, for an arbitrary
    $C_0$-semigroup
 $(e^{tA})$ on $H$, there exists $M \ge 1$ and $\omega \in \R$
   such that $\|
 e^{tA}\|_{{\cal L}(H)}
  \le M e^{\omega t}$, $t \ge 0.$ Moreover, in the sequel,
    $Lip(F)$ will  denote
   the
 Lipschitz constant of $F.$

\bpr \label{det} Let $F: H \to H$ be Lipschitz continuous and
bounded and $f \in L^p(0,T ; H)$, for some  $p >0$. Let $A:
D(A)\subset H \to H$ be the generator of a $C_0$-semigroup
$(e^{tA})$.

\hh (i) For any $x \in H$, the equation \beq \label{equ}
 y(t) = e^{tA} x + \int_{0}^t e^{(t-s) A } F(y(s) + f(s) ) ds
\eeq has a unique continuous solution $y: [0,T] \to H$.

\hh (ii)  There exists a constant $C= C (p, \omega, M, Lip(F), \|
F\|_0)>0$ such that for solutions $y$ and  $ z \in
 C([0,T]; H)$,   corresponding
 respectively to functions $f$, $g$  $ \in
    L^p(0,T ; H)$ and to the same $x \in H$,
  we have the estimates
 $$
 (a) \;\;\;
  \| y- z   \|_{C([0,T]; H)} \le C    \Big (\int_0^T | f(t) - g(t)|^p dt
  \Big)^{1/p}, \;\;\; p \ge 1;
$$ $$ (b) \;\;\;
\| y- z   \|_{C([0,T]; H)} \le C \int_0^T | f(t) - g(t)|^p
dt,\;\;\;  p \in (0,1).
$$
\epr

\bpf Assertion (i) follows easily by a fixed point argument. Let
us consider (ii). The proof of (ii) when $p \ge 1$ is an easy
 application of the Gronwall lemma. Thus we only  prove (b).

\vskip  1mm
  We consider a family  of equivalent norms $\|  \cdot \|_{\lambda}$ on
  the Banach space $E = C([0,T]; H)$, for $\lambda \ge 0$,
$$
  \| h \|_{\lambda} = \sup_{t \in [0,T]} e^{- \lambda t }
  |h(t)|,\;\;\; h \in E
$$
(for $\lambda =0$ we get the usual sup norm). For a fixed  $f
  \in L^p(0,T ; H) $, let us define the operator
   $K_f : E \to E$,
 $$
 (K_f  y) (t) = e^{tA} x + \int_{0}^t e^{(t-s) A } F(y(s) + f(s) )
 ds, \;\;\; y \in E, \;\; t \in [0,T].
$$
We find for $\lambda > \omega$,
    for any $y, z \in E$,
    \begin{align*}
 \| K_f y - K_f z   \|_{\lambda}
 \le C \sup_{t \in [0,T] } e^{- \lambda t }
\int_{0}^t e^{\omega (t-s) A } |y(s) - z(s)| ds
\\
 \le C \|y - z \|_{\lambda} \, \sup_{t \in [0,T] } \int_0^t e^{- (\lambda - \omega) (t-s) }
 ds  \,  \le \frac{ C}{\lambda - \omega }
\|y - z \|_{\lambda}.
\end{align*}
 Let us choose $\lambda_0$ large enough such that
 $c_0 = \frac{ C}{\lambda_0 - \omega } < 1$. We have
 \beq \label{fu}
\| K_f y - K_f z   \|_{\lambda_0} \le c_0
 \|y - z \|_{\lambda_0},\;\;\;\; y, z \in E.
  \eeq
Let now $f$ and $g \in L^p(0,T ; H)$. We get, for any $t
 \in  [0,T]$, $ y \in E$,
$$
|(K_f y) (t) -  (K_g y) (t) | \le M  \int_0^t e^{\omega (t-s)}
 | F(y(s) + f(s) ) - F(y(s) + g(s) )| ds.
$$
 Since $F$ is bounded and Lipschitz continuous, it is also
  H\"older continuous of order $p \in (0,1)$ and we find
$$
\| K_f y   -  K_g y  \|_{\lambda_0}
 \le  c M  e^{\omega T } \int_0^T
 | f(s)  -  g(s) |^p ds.
$$
 If we have solutions $y$ and $z$ corresponding to $f$ and $g$,
 then
 $y = K_f y $ and $z = K_g z$. We get
$$
 \| y - z\|_{\lambda_0 } =   \| K_f y - K_g z \|_{\lambda_0 }
$$
$$
 \le  \|  K_f y - K_f z \|_{\lambda_0 }  +
 \| K_f z - K_g z \|_{\lambda_0 } \le
  C_T \int_0^T
 | f(s)  -  g(s) |^p ds  + c_0   \|y - z  \|_{\lambda_0 }
 $$
and the assertion follows since $c_0 \in (0,1)$.

\epf

\begin{remark} { Clearly the previous result holds when $F$
 is only Lipschitz continuous and  $f \in L^p(0,T ; H)$ with $p \ge 1$.}
\end{remark}

\bth \label{exist} Assume Hypotheses \ref{basic} and \ref{dl}.
   Then there exists a unique  mild solution $(X_t^x)$
 to the equation \eqref{ab1}.
     Moreover $(X_t^x)$ is
   a Markov-Feller process.
 \eth
 \bpf  {\it Step 1. Existence and uniqueness}.
 Uniqueness follows by the Gronwall lemma. Let us prove
 existence.
  By using Proposition \ref{det} and Theorem \ref{bo}
   we find that, for any
   $x \in H$, there exists
 a continuous ${\cal F}_t$-adapted  process
   $(Y_t) =( Y_t^x)  $  with values in $H$ which solves  $\P$-a.s.
$$
  Y_t = e^{tA} x + \int_0^t e^{(t-s)A} F(Y_s + Z_A (s)) ds,\;\;
 \; t\ge 0.
 $$
 Let us define
$$
  X_t^x = Y_t^x +  Z_A (t),\;\;\; t \ge 0, \; x \in H.
 $$
  Since  $Z_A(t)$ is predictable it follows that
  $X = (X_t^x)$  is predictable as well.
Clearly $(X_t^x)$ is the unique mild solution.

\hh {\it Step 2. Markov property}.  The proof of the Markov
 property is quite involved. Indeed since our solution is not
 assumed to have  c\`adl\`ag trajectories,
 we have to
 proceed differently from \cite[Theorem 7.10]{DZ}.

For any measurable function $\psi : [0,T] \to H$, let $y(t)$
 be the unique continuous function with values in $H$
   which solves
 the equation
 $$
y(t) =  \int_0^t e^{(t-s)A} F(y(s) + \psi(s)) ds.
 $$
Set $y (t)= y(t, \psi )$, $t \in [0,T]$,  to stress the dependence
on $\psi$. We have
\begin{align*}
y (t+ h , \psi )= \int_0^{t+h} e^{(t+ h -s)A} F(y(s, \psi) +
\psi(s)) ds \\
= e^{hA} \int_0^{t} e^{(t -s)A} F(y(s, \psi) + \psi(s)) ds +
\int_t^{t+h} e^{(t+ h -s)A} F(y(s, \psi) + \psi(s)) ds
\\
= e^{hA} [y(t, \psi )] + \int_0^{h} e^{(h -s)A} F(y(t+ s, \psi) +
\psi(t+ s)) ds,\;\;\;\; t, t+h \in [0,T].
\end{align*}
Define a new function on $[0, T-t]$,
$$
v(\cdot, \psi  ) := y(t+ \cdot, \psi) - e^{ (\cdot )\, A} [y(t,
\psi )].
$$
We have
$$
v(h, \psi) = \int_0^{h} e^{(h -s)A} F( v(s, \psi) + e^{sA} [y(t,
\psi)] +   \psi(t+ s)) ds, \;\; h \in [0, T-t].
$$
By uniqueness, $ v(h, \psi) = y(h, e^{ (\cdot )\, A} [y(t, \psi )]
+ \psi (t+ \cdot)) $ and so we get $\text{for} \; t, t+h \in
[0,T],$
 \beq \label{f3}
 y \big( h, e^{ (\cdot )\, A} [y(t, \psi)] + \psi (t+ \cdot) \big)
 + e^{ h \, A} [y(t, \psi)]
= y(t+h , \psi).
 \eeq
Defining $u(t, \psi) = y(t, \psi ) + \psi(t)$, $t \in [0,T]$, we
find
\begin{align*}
& u(t+ h, \psi ) - \psi (t+h) = u \big(
 h, e^{(\cdot)A}[u(t, \psi) -
\psi (t) ] + \psi (t+ \cdot) \big)
\\
& - e^{hA}[u(t, \psi) - \psi (t) ] - \psi (t+h) + e^{hA}[u(t,
\psi) - \psi (t) ]
\end{align*}
and so we get
$$
u(t+ h, \psi ) = u (h, e^{(\cdot)A}[u(t, \psi) - \psi (t) ] + \psi
(t+ \cdot) )
$$
which is the same formula at  the end of  \cite[page 256]{DZ}. Now
the Markov property follows arguing as in \cite[page 257]{DZ}.

\vv {\it Step 3. Feller property.}   Fix $x, y \in H$. From the
estimate,
 \beq \label{r3}
|X_t^ x - X_t^y| \le \| e^{tA}\|_{{\cal L}(H)} \, |x- y| +
\int_0^t \| e^{(t-s)A} \|_{{\cal L}(H)} \, |F(X_s^x) - F(X_s^y)|
ds,
 \eeq
using the Lipschitz continuity of $F$ and the Gronwall lemma, we
find that, for any $T>0, $
   $|X_t^ x - X_t^y| \le M_T |x-y|$, $t
\in [0,T]$, $x, y \in H$, $\P$-a.s.. The Feller property (i.e.,
 the mapping $x \mapsto \E[f(X_t^x)]$ is continuous on $H$,
  for any $f\in C_b(H)$, $t \ge 0)$
follows easily. \epf

\subsection{Irreducibility}

We establish now irreducibility of the solutions to \eqref{ab1}.
\bth \label{irre} Assume Hypotheses \ref{basic} and \ref{dl}.
Then, for any $x \in H$, the mild solution $X=(X_t^x)$
 to the equation \eqref{ab1} is irreducible.
\eth \bpf
 Fix $x \in H$, $T>0,$ and denote by  $X = (X_t)$ the
solution to \eqref{ab1} starting from $x$. Set
 $$
 Y_t = X_t - Z_A(t),\;\; t \in [0,T],
 $$
where \beq\label{a4} \left\{
\bal dZ_A(t)  & = A Z_A(t) dt  + dZ_t, \\
      Z_A(0) & = 0, \;\;\;  t \ge 0.
\eal \right. \eeq Note that
$$
 Y_t = e^{tA} x + \int_0^t e^{(t-s)A} F(Y_s + Z_A (s)) ds.
$$
 Let $z^{u}$ and $y^{u,x}$ be
  the solutions,  driven by a control function $u\in L^2(0,T ; H)$,
   of the following
 control systems, respectively,
  \beq\label{a3} \left\{
\bal \frac{dz}{dt}  & = Az(t) \,  +  u(t), \\
      z(0) & = 0, \;\;\;  t \in [0,T],
\eal \right.    \;\;\;\;\; \;\;\; \left\{
\bal \frac{dy}{dt}  & = Ay(t) + F(y(t))\, +  u(t), \\
      y(0) & = x \in H, \;\;\;  t \in [0,T].
\eal \right.
 \eeq   Thus
\begin{equation}
z^{u}(t) =  \int_{0}^{t} e^{(t-s)A} u(s)ds ,\,\,\,t\in [0,T],
\end{equation}
and  $y^{u,x}$ is the solution of the following integral equation
$$
 y(t) = e^{tA} x + \int_0^t e^{(t-s)A} F(y(s))
 ds + z^{u} (t),\,\,\,t\in [0,T].
$$
By Theorem 7.4.2 of \cite{DZ1} we know that the second system in
(\ref{a3}) is approximately controllable at time $T>0$ in the
sense that, for any $x, a \in H$ and for any  $\epsilon >0$, there
exists  a control function $u\in L^2(0,T ; H) $ such that $|
y^{u,x} (T) -a| < \epsilon $.

Let
$$
 \bar y (t) = y^{u,x} (t) -  z^{u} (t),\;\;\; t \in [0,T].
$$
Note that
$$
 \bar y(t) = e^{tA} x + \int_0^t e^{(t-s)A} F(\bar y(s) +
  z^{u} (s)) ds.
$$
 Take $p < \alpha$ with $p \in (0,1)$. By estimate (b) in
 Proposition \ref{det} we get, $\P$-a.s.,
$$
 \sup_{t \in [0,T]} | Y_t - \bar y(t)   |  \le C \int_0^T | Z_A(t) -
z^{u}(t)|^p dt.
$$
 and so $| Y_T - \bar y(T)   |  \le C \int_0^T | Z_A(t) -
z^{u}(t)|^p dt$ or, equivalently,
$$
|X_T - Z_A(T) - y^{u,x} (T) +  z^{u} (T)  | \le C \int_0^T |
Z_A(t) - z^{u}(t)|^p dt.
 $$
We write, for any $a \in H$,
\begin{align*}
|X_T - a | \le |X_T - Z_A(T) - y^{u,x} (T) +  z^{u} (T)  |
 +
| Z_A(T) + y^{u,x} (T) -  z^{u} (T) -a |
\\
\le C \int_0^T | Z_A(t) - z^{u}(t)|^p dt + | y^{u,x} (T) -a| +
  |Z_A(T) -  z^{u} (T)| = I_1 + I_2 + I_3.
\end{align*}
For a given $\epsilon >0$, let us fix a control function $u$
 such that $I_2 = | y^{u,x} (T) -a| < \epsilon /3$.
  Using Proposition \ref{bo}, we get with positive probability
   that $I_1 < \epsilon /3$ and $I_3 < \epsilon /3$.
 The result follows.
\epf

\subsection{Strong Feller property }

  Let $(P_t) $ be the Markov semigroup associated to
   $X = (X_t^x)$, i.e.
 $P_t : B_b (H) \to B_b (H)$,
 \beq \label{nn}
 P_t f (x) = \E [f(X_t^x)],\;\;\; x\in H,\;\; f \in B_b(H), \;\;
 t\ge 0.
 \eeq
 To show the strong Feller property of $(P_t)$,
  we will assume
 Hypotheses \ref{basic}, \ref{dl} and
\begin{hypothesis} \label {basic3}   Assume that $\alpha \in (1,2)$ and that
 there exists  $\gamma \in [1/\alpha,1) $ and
  $C_1 >0$ such that
 \begin{equation} \label{pu1}
 \beta_n \ge C_1 \,
 \gamma_n ^{( \frac{1}{\alpha}) - \gamma},\;\;\; n \ge 1.
\end{equation}

\end{hypothesis}
 \noindent   This assumption is stronger than Hypothesis
 \ref{basic1}. Indeed, assuming
 (i) in  Hypothesis \ref{basic},  Hypothesis \ref{basic3} holds
  if and only
 there exist $ \hat c >0$,  $\gamma \in [1/\alpha,1)$,   such that
 \begin{equation} \label{pu}
 \sup_{n \ge 1} \frac{
\gamma_n ^{1/ \alpha} \alpha^{1/\alpha}} { \beta_n^{}
\, (e^{\alpha \gamma_n t} -1)^{1/\alpha}}
 \le \frac{\hat c} {t^{\gamma}} ,  \;\; \;\;\;  t >0
\end{equation}
and according to \eqref{ritorno} this implies \eqref{e4}.
 To see the previous equivalence  let us first assume \eqref{pu}. By choosing $t =
  \frac{1}{\gamma_n}$, we find \eqref{pu1}, with
  a suitable constant $C_1$.
 Viceversa, assume \eqref{pu1} with $\gamma \in [1/\alpha, 1)$.  Setting $d=
 \sup_{r>0} \frac{r^{\gamma}}{(e^r -1)^{1/\alpha}}$ we find
$$
\sup_{n \ge 1} \frac{
t^{\gamma}\, \gamma_n ^{1/ \alpha} \alpha^{1/\alpha}} { \beta_n^{}
\, (e^{\alpha \gamma_n t} -1)^{1/\alpha}}
 \le \alpha^{1/\alpha - \gamma} \, d \, \frac{\gamma_n^{1/\alpha - \gamma}}{\beta_n} \le
 \alpha^{1/\alpha - \gamma} \, \frac{d}{C_1}.
$$

 Before stating our  theorem on the strong
 Feller property,
 we discuss a motivating  example.

\begin{example}\label{E2} {\em Consider
  the following non-linear version of the  stochastic  heat
  equation  on $D= [0, \pi]^d$ treated in Example \ref{E1}:
 \beq \label{heat12}
  \left\{
\bal dX(t, \xi)  &  = \triangle X(t, \xi) \,dt  +
 f(X(t, \xi))dt \, + \,dZ(t, \xi), \;\;\; t>0,\\
      X(0, \xi) & = x(\xi), \;\;\; \xi \in D,
\\ X(t, \xi) &=0, \;\; t>0, \;\;\; \xi \in \partial D,
\eal \right. \eeq
 where $f: \R \to \R$ is a bounded and Lipschitz continuous
 function and
  $Z= (Z_t)$ is a  cylindrical $\alpha$-stable process
 $$
 Z_t = \sum_{j = (n_1, \ldots , n_d) \in \N^d}
  \beta_j Z_t^j e_j,
$$
 $\alpha \in (1,2)$,  where  $(e_j)$ is the  basis of  eigenfunctions of the Laplacian $\Delta$ in $H = L^2(D)$ (with  Dirichlet boundary conditions).
Thus if
 \begin{equation} \label{fy}
\bal & \sum_{j  \in \N^d} \beta_j^{\alpha} < + \infty \;\;\;
  \mbox{and}
 \;\;\; \beta_{(n_1, \ldots , n_d)}
  \ge c (n_1^2+  \ldots  + n_d^2)^{ {1/ \alpha - \gamma}},
   \;\; (n_1, \ldots , n_d) \in \N^d,
\eal
\end{equation} for some constants  $c >0$ and  $\gamma \in [1/\alpha,1)$,
 then, by Theorems \ref{irre} and \ref{main}, the solution
  to \eqref{heat12} is irreducible and strong Feller. By  the
   Doob theorem (see \cite[Theorem 4.2.1]{DZ1}) if there
    exists an invariant measure for \eqref{heat12} then it must be
    unique.
In particular,   if
 $ \beta_{(n_1, \ldots , n_d)}
  =  (n_1^2+  \ldots  + n_d^2)^{\, 1/ \alpha -\gamma},
  $ $ \;\; (n_1, \ldots , n_d)$ $ \in \N^d,$ then \eqref{fy}
 is equivalent  to
$$
\sum_{ (n_1, \ldots , n_d)  \in \N^d}  (n_1^2+  \ldots  +
n_d^2)^{1- \alpha \gamma } < + \infty.
$$
 This   holds if and
 only if $\alpha \gamma > \frac{d}{2} + 1$. Thus if $d=1$,
 one requires that $\alpha \gamma > \frac{3}{2}$.}
 \qed
\end{example}

\bre {\em The above  example shows that the strong Feller
 property holds in rather special situation. It seems thus of great
 interest to develop the concept of asymptotic strong Feller
 property in the case of SPDEs with L\'evy noise
 (compare with \cite{HM1}, \cite{HM2} and \cite{HM3}).
 } \ere

 \bth \label{main}  Assume that Hypotheses \ref{basic},
\ref{dl}
 and \ref{basic3} hold.
 Then, for any $t>0,$ the transition semigroup $(P_t)$
  (see \eqref{nn}) maps Borel and bounded
 functions
 into Lipschiz continuous functions.  Moreover, there exists
  $ \tilde C = \tilde C (\gamma, c_{\alpha}, C_1,
\| F\|_0
 )>0$, such that,
  for any $x,y \in H$,
   we have
 \beq \label{fw1}
 |P_t f(x) - P_t f(y)| \le  \tilde C   \| f\|_0
 \,  \frac{1}{t^{\gamma} \wedge  1}  |x-y|,\;\;\;\;\;
   t >0,\;\;\; f \in B_b (H).
  \eeq
\eth

 Recall that  $t^{\gamma} \wedge 1 = {\min (t^{\gamma} , 1)}$.
  Note that it is enough to prove
  estimate \eqref{fw1}, for any $t \in (0, 1]$. Indeed,
   when   $t >1$,
  we can replace  $|P_t f(x) - P_t f(y)|$
   in \eqref{fw1}  with
  $|P_{1} (P_{t-1} f) (x) - P_{1} (P_{t-1} f)
 (y)|$.

\smallskip To prove the result we first investigate generalised
   solutions to the Kolmogorov equation associated to $(P_t)$ (or to
   $(X_t^x)$) as in
   \cite[Section 9.4.2]{DZ}.

 Note that the   generator ${\cal A}_0$ of $(P_t)$ is formally given by
  \beq \label{kol}
 {\cal A}_0 f (x) = \lan Ax + F(x), D f(x)\ran  + \sum_{n \ge 1}
  \beta_n^{\alpha} \int_{\R} (f (x+ e_n z) - f(e_n z))
\frac{1}{|z|^{1+
  \alpha}}dz,
 \eeq
 for regular and cylindrical functions $f: H \to \R$.
  The associated
 Kolmogorov equation is
\begin{equation} \label{eco} \left\{
\begin{aligned}  \partial_t u(t,x) & = {\cal A}_0 u (t, \cdot)(x), \;\;\; t>0,\;
x \in H,\\
      u(0, x) & = f(x),\;\; x \in H.
\end{aligned}  \right. \end{equation}
 Let us fix $T>0$ and consider the space
 $$
 \Lambda (0,T) = \{ u \in C(]0,T]; C_b^1 (H))\;\; :\;\;
  \sup_{t \in ]0,T]} t^{\gamma} \|u(t, \cdot) \|_1 < \infty
   \},
 $$
 where $ \|u(t, \cdot) \|_1 =   \| u (t, \cdot)\|_0
  +   \|D_x u(t, \cdot) \|_0 $ and
  $\gamma \in (0,1)$ is fixed in   Hypothesis
 \ref{basic3}. According to \cite{DZ}
 a {\it mild
 solution}  to the Kolmogorov equation \eqref{eco}
  (on $[0,T]$
   with initial datum $f \in B_b (H))$ is a function $u \in
 \Lambda (0,T)$ such that
 \beq \label{km}
 u(t,x) = R_t f (x) + \int_0^t R_{t-s}
\big (  \lan F(\cdot), D u(s, \cdot ) \ran\big) (x)\, ds , \;\; t
\in [0,T],\; x \in H, \eeq
 where $D = D_x$ and $(R_t)$ is the transition semigroup
 determined by
   the
 linear equation \eqref{ou}. To stress the dependence on $f$, we
 will also write
 $$
u = u(t,x) = u^{f}(t,x),\;\;\; t \in [0,T],\;\; x \in H.
$$
 Note that using Theorem \ref{grad} and Hypothesis \ref{basic3}
(see also \eqref{pu}), we get, for any  $f \in B_b (H)$,
 \beq \label{rt}
 \| DR_t f \|_0 \le \frac{C_0}
  {t^{\gamma}} \| f\|_0,\;\;\; t >0,
  \;\; \mbox{where} \;\; C_0 = 8 c_{\alpha} \, \hat c. \eeq
 Thanks to \eqref{rt}, we can adapt the proof of  \cite[Theorem 9.24]{DZ}
 and obtain  that the mapping
 $S : \Lambda(0,T) \to \Lambda (0,T)$,
 \beq \label{sss}
S (u) (t,x) = R_t f (x) + \int_0^t  R_{t-s} \big (  \lan F(\cdot),
D u(s, \cdot ) \ran\big)(x) ds,\;\; u \in \Lambda(0,T),
 \eeq
is a contraction for $T$ small enough. Therefore, we obtain

\bpr \label{univ} For any $f \in B_b(H)$, $T>0$, there exists a
 unique mild solution $u = u^f$ to \eqref{eco}.
  Moreover, for any $t \ge 0$, we may define:
$$
\tilde P_t  f (\cdot) := u^{f}(t, \cdot), \;\; \;\;
 f \in B_b
(H).$$
 It turns out that  $(\tilde P_t)$ is a semigroup  of bounded linear
 operators on $B_b (H)$.
\epr \noindent

In the proof  of the next lemma,  we will use the following
 Gronwall type lemma. Let $a, b, {\gamma} $ be non-negative constants,
with $ {\gamma} <1$. Let $T>0$. For any integrable function $v:
[0,T] \to \R$,
 \beq \label{gron}
 0 \le v (t) \le  a t^{-{\gamma}} + b\int_0^t (t-s)^{-{\gamma}} v(s) ds, \;\;
 t \in [0,T[ \; \; a.e.,\;\;\; \mbox{implies} \;\; \;
  v(t) \le a M t^{-{\gamma}},
\eeq
 $t \in [0,T[ $, a.e.. (where $M =M (b,{\gamma},T)   1 + b \, k_{{\gamma}} \, T^{1-{\gamma}}$).

\ble \label{apri}  For any $T>0$, there exists $c= c (\gamma,
c_{\alpha}, C_1, \| F\|_0,
 T)>0$ such
that, for any $f \in B_b (H)$, $t \in ]0,T]$,
$$
  \| D \tilde P_t f  \|_0 \le \frac{c}{t^{\gamma}} \| f\|_0.
 $$
\ele
 \noindent {\it Proof.} We have
$$
 Du (t,x) = DR_t f (x) + \int_0^t D R_{t-s}
\big (  \lan F(\cdot), D u(s, \cdot ) \ran\big) (x)\, ds , \;\;\; x
\in H.
$$
 By using \eqref{rt} and the previous Gronwall lemma, we
 get
$$
\|  Du (t, \cdot ) \|_0 \le  \frac{C_0 M } {t^{\gamma}} \|  f
\|_0, \;\; t \in ]0, T],\;\;\; M = M (\gamma, c_{\alpha}, \hat c,
\| F\|_0, T)>0. \qed
$$
 \vv \noindent {\bf Galerkin's approximation.} To show the
regularizing effect of $(P_t)$, according to \cite[Theorem
9.27]{DZ},
   it would
be enough  to prove   that
  $(P_t)$
 and $(\tilde P_t)$ coincide. However
 the proof of \cite[Theorem 9.27]{DZ} is not complete and we
  are unable to fill the gap in our  situation.

   We therefore  resort to Galerkin's
  approximations and we will only  identify suitable
 finite-dimensional   semigroups which approximate
   $(P_t)$
 and $(\tilde P_t)$ respectively.

  \vv Let us consider orthogonal projections $\pi_n : H \to H_n$,
   $n \in \N,$ where $H_n$  is the subspace of $H$
   generated by $\{ e_1,
\ldots, e_n \}$. For any $n \in \N$, $x \in H$, define the
$H_n$-valued process  $(Y^n_t) =(Y^n_t(x))$ as the unique mild
solution to
 \beq \label{fiii}
  \bal
 Y_t^n  = e^{tA_n} x + \int_0^t e^{(t-s)A_n} (\pi_n \circ F \circ \pi_n)
  ( Y_s^n) ds +  Z_{A_n} (t),\;\;\;
 \eal
  \eeq
where $A_n = \pi_n \circ A $. Let  $F_n = \pi_n \circ F \circ
\pi_n $. Note that, for any $n \in \N$, it holds:
 \beq \label{fr}
 \| F_n \|_0  \le \| F\|_0,\;\;\;\; Lip (F_n ) \le Lip (F),
\eeq
 where $Lip (F_n)$ denotes the Lipschitz constant of $F_n$.

\vv Consider the mild solution $u_n$ to the Kolmogorov equation
 corresponding to $Y_t^n$, i.e.,
 \beq \label{f4}
 \bal
u_n (t,x) = u_n^f (t,x)
 = R_t^n f (x) + \int_0^t  R_{t-s}^n \big (
\lan F_n(\cdot), D u_n(s, \cdot ) \ran\big)(x) ds , \;\;\; x \in
H,
\\
 \mbox{where} \;\; R_t^n f(x) = \E [f(e^{tA_n} x + \pi_n Z_A(t))]
  = \int_{H} f(e^{tA_n}x + \pi_n y  ) \mu_t^0 (dy).
\eal
 \eeq
Define the following
 two approximating semigroups on $B_b (H)$ (see
\eqref{fiii} and \eqref{f4}):
  \beq
 P_t^n f (x) = \E[f(Y_t^n(x))],\;\;\; \tilde P_t^n f(x) = u_n^f (t,x),
 \;\; f \in B_b(H),
 \eeq
\ble \label{nnn} For any function $f \in B_b (H)$, $n \in \N$, we
have
$$
P_t^n f = \tilde P_t ^n f ,\;\;\; t \ge 0.
$$
\ele \bpf We fix $n \in \N $. It is enough to prove the assertion
for any cylindrical function $f \in B_b(H)$, which depends only on
the first $n$-coordinates. Identifying $(P_t^n)$  and  $(\ti P_t
^n )$ with the corresponding semigroups acting on $B_b(\R^n)$, we
have to check that
 \beq \label{tesi}
P_t^n f = \ti P_t ^n f ,\;\;\;f \in B_b(\R^n),\;\;
 t \ge 0.
  \eeq
To this purpose (identifying $F_n$ with the corresponding
 Lipschitz continuous function from $\R^n$ into $\R^n$)
first note  that $(P_t^n) $ is a    strongly continuous semigroup
of positive contractions on
 $C_0(\R^n)$
 (see \cite[Section 6.7]{A}). Here $C_0(\R^n)$ denotes the space of all real continuous
  functions on $\R^n$ vanishing at infinity.

\smallskip Let us consider now $(\ti P_t^n)$. We start to show that
$
 \ti P_t^n (C_0(\R^n)) \subset C_0(\R^n)$, $t \ge 0$.

Fix $T>0$ and let $f \in C_0(\R^n)$ and $t \in ]0,T]$; we will use
an
 inductive argument to prove that $\ti P_t f \in C_0(\R^n)$.
 By \eqref{sss}, we know that
 $$
 \ti P_t^n f = \lim_{m \to \infty} S^m (0)
  = \lim_{m \to \infty} (S \circ \ldots \circ S) (0) \;\;\;\; \mbox{in} \;\;
 $$ $$  \Lambda
  (0,T) =  \{ u \in C(]0,T]; C_b^1 (\R^n)) :\;\;
  \sup_{t \in ]0,T]} t^{\gamma} \|u(t, \cdot) \|_1 < \infty
   \}.
$$
 We prove that, for any $m \in \N$, $S^m (0) (t, \cdot )$ and
  $D_x S^m (0) (t, \cdot ) \in C_0 (\R^n)$. We have (for $m =1$)
 $S^1 (0) (t, \cdot )(x)$ $ = R_t f (x)$, and so
 \beq \label{der}
  \bal
& D_x S^1 (0) (t, \cdot )(x) = DR_t^n f (x)
 = \int_{\R^n} f(e^{tA_n} x +  y ) \, U_n(y,t) \; \mu_t^n (dy),
 \;\; x \in \R^n, \; \mbox{where}
\\
& \mu^n_t \;\; \text{has  density}
 \;\; \prod_{k=1}^n
  p_{\alpha} \Big(\frac{ y_k }{c_k (t)}
 \Big) \frac{1}{c_k (t) } \;\;
 \\ & \text{ and }
  \;\; U_n(y,t) = \sum_{k = 1}^n \,
  \frac{p'_{\alpha}(\frac{y_k} {c_k (t)})}
  {p_{\alpha}(\frac{y_k} {c_k (t)})} \,
\frac{e^{- \gamma_k t }}{c_k (t)} \, e_k
  \; \in L^2 (\mu^n_t; \R^n).
\eal
 \eeq
It follows easily that  $S^1 (0) (t, \cdot )$ and
  $D_x S^1 (0) (t, \cdot ) \in C_0 (\R^n)$. Assume that the
  assertion holds for an arbitrary  $m \in \N$. Since
 $$
 \bal
 & S^{m+1} (0) (t, \cdot )(x)
  = R_t^n f (x) + \int_0^t  R_{t-s}^n \big (  \lan F_n
   (\cdot), D  S^{m}(0)(s,
\cdot ) \ran\big)(x) \, ds,
 \\
& D_x S^{m+1} (0) (t, \cdot )(x) =   D R_t^n f (x) +
  \int_0^t  D R_{t-s}^n \big (  \lan F_n(\cdot), D  S^{m}(0)(s,
\cdot ) \ran\big)(x)\, ds  = D R_t^n f (x)
\\ & + \int_0^t ds
  \int_{\R^n}  \lan F_n(e^{(t-s) A_n} x +  y ), D  S^{m}(0)(s,
  e^{(t-s) A_n} x +  y) \ran \,
U_n(y,t-s) \; \mu_{t-s}^n (dy),
 \eal
$$
$x \in \R^n,$ we have easily that the assertion
 holds also for $m +1$.

\hh Using  Lemma \ref{apri}, we get that $(\ti P_t^n)$ is a
 strongly continuous
  semigroup of  bounded linear operators on $C_0(\R^n)$.

\vv We will  prove \eqref{tesi} when $f \in C_0(\R^n)$. Indeed, by
a standard argument (see \cite[Chapter 4]{EK}) this is enough to
get \eqref{tesi}.

\vv By Ito formula $D_0 = C^2_0 (\R^n) = \{ f \in C_0 (\R^n) \,\;
: \; \, Df \; and \;\; D^2 f \in C_0(\R^n) \}$
 is invariant for $(P_t^n)$ (compare with \cite[Theorem 6.7.4]{A}).
 Moreover,   $D_0 \subset $dom $({\cal A}_n)$, where ${\cal A}_n$ is
the generator of $( P_t^n)$. By a well known result, $D_0$ is a
core for $( P_t^n)$ (see \cite[page 52]{EN}). Note that
$$
 {\cal A}_n f (x)= \lan A_n x + F_n(x), D f(x)\ran  + \sum_{k = 1}^n
  \beta_k^{\alpha} \int_{\R} (f (x+ e_k z) - f(e_k z))
\frac{1}{|z|^{1+
  \alpha}}dz, \;\; f \in D_0.
$$
 Let us consider $(\ti P_t^n)$. If $f \in D_0$,
  we can solve (by the contraction principle)
$$
u (t,x) = R_t^n f (x) + \int_0^t  R_{t-s}^n \big (  \lan
F_n(\cdot), D u(s, \cdot ) \ran\big)(x)ds , \;\;\; x \in \R^n,
$$
in the space $C ([0,T]; C^2_0(\R^n))$ and get that $D_0$
 is also invariant for $(\ti P_t)$.
 A straightforward calculation, shows that $D_0
  \subset $dom $(\ti {\cal A}_n)$, where $\ti {\cal A}_n$ is
 the generator of $( \ti P_t^n)$. Thus $D_0$ is
a core also for $(\ti P_t^n)$. Moreover,    $\ti {\cal A}_n$
 coincides
 with $ {\cal A}_n$ on $D_0$.
 It follows that $(P_t^n) $ and $(\ti P_t^n)$ coincide on
$C_0(\R^n)$
and this
finishes the proof.
 \epf

 \newpage \noindent \textbf{\emph  Proof of Theorem \ref{main}.}

 \vskip 1mm \noindent  {\it I Step.} Using  \eqref{f4}, Lemmas
  \ref{apri} and
  \ref{nnn} and  the semigroup property, there exists
  $\tilde C = \tilde C (\gamma, c_{\alpha},  C_1,
\| F\|_0)>0$ such that,  for any $f \in C_b (H)$,
  \beq \label{gt}
 \bal
& |u_n (t,x) - u_n(t,y)| = |P_t^n f(x) - P_t^n f(y)|  \le |R_t^n f
(x) - R_t^n f (y)|  \\ & +
 \int_0^t | R_{t-s}^n \big (
\lan F_n(\cdot), D u_n(s, \cdot ) \ran\big)(x)
 - R_{t-s}^n \big (
\lan F_n(\cdot), D u_n(s, \cdot ) \ran\big)(y) | ds
\\ & \le
  \tilde C  \| f\|_0
  \frac{1}{t^{\gamma} \wedge 1}
   |x-y|,\;\; x, y \in H, \;\; n
  \in \N,\; t>0.
\eal
 \eeq
{\it II Step. } For any $f : H \to \R$ which is
 continuous and bounded, we have:
$$
 \lim_{n \to \infty } P_t^n f (x) = P_t f(x),\;\; x \in H,\;\; t
 \ge 0.
$$
 Recall that  $P_t^n f(x) = \E [f (Y_t^n(x) )]$
  (see
\eqref{fiii}). The assertion will  follow by proving that
 \beq \label{ce}
 \lim _{n \to \infty} Y_t^n(x) = X_t^x,\;\;\; x \in H,\;\; t \ge 0,
 \;\; \P-a.s..
\eeq To show \eqref{ce}, we fix $x \in H$ and write
 \begin{align*}
  X_t^{  x}  = e^{tA} x + \int_0^t e^{(t-s)A}
   F_n  ( X_s^{ x}) ds +  Z_{A} (t) +  f_n (t), \;\;\; \mbox{where}
\\
 f_n (t) = \int_0^t e^{(t-s)A}
   [F(X_s^{ x} ) - F_n  ( X_s^{ x})]  ds,\;\;\; t \ge 0, \; n \in
   \N.
\end{align*}
(see the notation in \eqref{fiii}). Defining
 $U_t^n = X_t^{ x} - Z_A(t) - f_n (t)$, we have:
$$
U_t^n =  e^{tA} x + \int_0^t e^{(t-s)A}
   F_n  ( U_s + Z_A(s) + f_n (s) ) ds.
$$
Note that
\begin{align*}
 Y_t^n(x)  = e^{tA} x + \int_0^t e^{(t-s)A} F_n
  ( Y_s^n(x)) ds +  Z_{A_n} (t)  + g_n (t),\;\;\; \mbox{where}
\\
 g_n (t) = e^{tA_n} x - e^{tA} x,\;\;\; t \ge 0, \;\; n \in \N.
\end{align*}
Introducing $V_t^n = Y_t^n(x) - Z_{A_n}(t) - g_n (t)$, we find
$$
 V_t^n =  e^{tA} x + \int_0^t e^{(t-s)A}
   F_n  ( V_s^n + Z_{A_n}(s) + g_n (s) ) ds.
$$
Now to estimate $|U_t^n - V_t^n|$, we use (b) in Proposition
\ref{det}. Let us choose
 $p \in (0, \alpha) \cap (0,1)$ and fix any $T>0$. We have
\beq \label{r5}
 \sup_{t \in [0,T]} |U_t^n - V_t^n| \le
  C \int_0^T |Z_{A}(s) + f_n (s) -  Z_{A_n}(s) - g_n (s) |^p ds.
\eeq
Note that $ Z_{A_n}(t) = \pi_n Z_A(t)$, $n \in \N$. Moreover,
it is easy to see that
$$
\lim_{n \to \infty} |f_n (t)| + |g_n (t)| =0,
$$
for any $t \ge 0$. Applying the dominated convergence theorem in
 \eqref{r5}, we infer
$$
 \lim_{n \to \infty} \sup_{t \in [0,T]} |U_t^n - V_t^n| =0.
$$
Using the inequality
\begin{align*}
& |Y_t^n(x) - X_t^x| \le | Y_t^n(x) - Z_{A_n}(t) - g_n (t) -
 X_t^{ x} + Z_A(t) + f_n (t)| \\
  & + |Z_{A_n}(t) + g_n (t)
 - Z_A(t) - f_n (t)  |,
\end{align*}
$t \ge 0$, $n \in \N$, and passing to the limit as $n
\to \infty$, we get assertion \eqref{ce}.

\hh {\it III Step. } By  the previous steps we know that, for any
 $f \in C_b (H)$,
$$
|P_t f (x) - P_t f (y) |  \le
  \tilde C  \| f\|_0
  \frac{1}{t^{\gamma} \wedge 1}
   |x-y|,\;\; x, y \in H, \;\;  t>0.
$$
 Now we get the assertion, using that
$$
Var  \big [ p_t (x, \cdot) -  (p_t (y, \cdot) \big ]
 = \sup_{f \in C_b (H),
 \; \| f\|_0 \le 1} |P_t f (x) - P_t f (y) |,
$$
for any $t>0,$ $x, y \in H$, where
  $p_t(x, \cdot)$ denotes the kernel of $P_t$ and {\it Var}
  the
 total variation (see the proof of \cite[Theorem 9.28]{DZ} or
\cite[Lemma 7.1.5]{DZ1}).  \qed

{{\vskip 3mm \noindent}} {\bf Acknowledgment } The authors wish to
 thank professor Kwapien for providing the
  proof of part (ii) in  Theorem \ref{kw},
 which significantly improved  our
 initial result.

\end{document}